\documentclass[a4paper, 12pt, final]{amsart}

\title[Algebraically bounded derivations II]{Generic derivations on algebraically
  bounded structures II\\ Model theoretical properties}

\date{24/09/2024}

\usepackage[T1]{fontenc}
\usepackage{lmodern}
\usepackage{microtype}

\author[A. Fornasiero]{Antongiulio Fornasiero}
\address{Universit\`a di Firenze}  
\email{antongiulio.fornasiero@gmail.com} 
\urladdr{https://sites.google.com/site/antongiuliofornasiero/}

\author[G. Terzo]{Giuseppina Terzo}
\address{Universit\`a degli Studi di Napoli ``Federico II''}
\email{giuseppina.terzo@unina.it}

\newcommand*{\mykeywords}{Derivation, algebraically bounded, model completion,
  imaginaries, simplicity, open core}

\usepackage[svgnames]{xcolor}
\usepackage[pdftex, linktocpage=true, 
pdfborder={0 0 0}, plainpages=false,%
pdfpagelabels, pdfdisplaydoctitle=true,%
pdfstartview={XYZ null null null},%
colorlinks=true,
citecolor=ForestGreen,
unicode=true,
pdfusetitle
]{hyperref}
\hypersetup{%
  pdfauthor={A. Fornasiero, G. Terzo},
  pdfkeywords={\mykeywords},
  pdfsubject={Derivations on an algebraically bounded structure,
model completion of the corresponding theory}
}

\usepackage{amsmath, amsthm, amssymb}

\newtheorem*{thmA}{Theorem A}

\usepackage{aliascnt}
\let\oldtheorem\newtheorem
\RenewDocumentCommand{\newtheorem}{s m o m O{}}{%
\IfBooleanTF{#1}%
{\oldtheorem{#2}{#4}}%
{\IfNoValueTF{#3}{\oldtheorem{#2}{#4}[#5]}%
{\newaliascnt{#2}{#3}%
\oldtheorem{#2}[#2]{#4}%
\aliascntresetthe{#2}}}}

\usepackage{mathtools}
\usepackage{xspace}

\usepackage{enumitem}

\setlist[itemize,1]{
  wide,
}

\setlist[enumerate,1]{%
  wide,
  label = \arabic*), 
         font = \upshape%
       }

\SetEnumitemKey{Romanenum}{
  label = \Roman*), 
  font = \upshape,
}

\SetEnumitemKey{romanenum}{
  label = (\roman*), 
  font = \upshape,
}

\SetEnumitemKey{romanenump}{
  label = (\roman*'), 
  font = \upshape,
}

\newlist{customenum}{enumerate}{1}
\makeatletter
\def\ctext#1{\expandafter\@ctext\csname c@#1\endcsname}
\def\@ctext#1{\ifcase#1\or Deep\or Wide\or PP\or
Fourth\or Fifth\or Sixth\fi}
\makeatother
\AddEnumerateCounter{\ctext}{\@ctext}{Deep}
\setlist[customenum]{label=(\ctext*), 
  ref=(\texttt{\ctext*}),
  font= \texttt
}

\usepackage[english]{babel}
\usepackage[babel]{csquotes}
\usepackage[backend=biber, style=alphabetic, sorting=nyt,
 backref = true, backrefstyle = three,
 isbn = false, 
 ]{biblatex}
\AtEveryBibitem{%
  \clearlist{language}
} 

\renewbibmacro*{doi+eprint+url}{%
  \iftoggle{bbx:doi}
    {\printfield{doi}}
    {}%
  \newunit\newblock
  \iftoggle{bbx:eprint}
    {\usebibmacro{eprint}}
    {}%
  \newunit\newblock
  \iffieldundef{doi}
    {\iftoggle{bbx:url}%
      {\usebibmacro{url}}
      {}}
    {}%
}

\usepackage[
capitalize]{cleveref}

\newcommand{\Thms}{Thm.~}

\newcommand{\Defs}{Def.~}
\newcommand{\Props}{Prop.~}
\newcommand{\Cors}{Cor.~}
\newcommand{\Rems}{Rem.~}
\newcommand{\Conjs}{Conj.~}

\DeclareMathOperator{\RCF}{RCF}
\DeclareMathOperator{\ACFz}{ACF_0}
\DeclareMathOperator{\DCF}{DCF}
\DeclareMathOperator{\ldim}{loc-dim}

\DeclareMathOperator{\sing}{Sing}
\DeclareMathOperator{\reg}{Reg}
\DeclareMathOperator{\td}{tr.deg.}

\DeclareMathOperator{\MoR}{MR}

\newcommand{\ordefinable}{$\mathord{\vee}$\hyph definable\xspace}
\newcommand{\realcanon}{real\hyph canonical\xspace}
\newcommand{\eqcanon}{canonical\xspace}
\newcommand{\realstrict}{real\hyph strict\xspace}
\newcommand{\eqstrict}{strict\xspace}

\DeclareMathSymbol{\mlq}{\mathord}{operators}{``}
\DeclareMathSymbol{\mrq}{\mathord}{operators}{`'}

\newcommand{\IACS}{IACS\xspace}
\newcommand{\GEI}{GEI\xspace}

\newcommand{\Gdelta}{\mathcal G_{\delta}}
\DeclareMathOperator{\Diag}{Diag}

\newcommand{\K}{\mathbb K}

\newcommand{\Kt}{\tilde K}

\newcommand{\Kdb}{\tuple{\K, \bar \delta}}

\newcommand{\Td}{T^{\delta}}
\newcommand{\Tdg}{T^{\delta}_g}

\newcommand{\Tdgpp}{{T^{\delta}_{\mathrm{wide}}}}

\newcommand{\Tdb}{T^{\bar \delta}}
\newcommand{\Tdbg}{T^{\bar \delta}_g}
\newcommand{\Tdn}{T^{\bar \delta,nc}}
\newcommand{\Tdng}{T^{\bar \delta,nc}_g}

\newcommand{\Tds}{T^{\bar \delta,?}}
\newcommand{\Tdsg}{T^{\bar \delta,?}_g}
\newcommand{\cfieldK}{\mathcal C_{\deltabar}}
\DeclareMathOperator{\DerK}{Der_{K}}
\newcommand{\taud}{\tau_d}
\newcommand{\tauFO}{\tau_{FO}}
\newcommand{\Generic}{\mathbb G}

\newcommand{\Ldb}{L^{\bar \delta}}
\newcommand{\deltab}{\bar\delta}
\newcommand{\epsb}{\bar\eps}
\newcommand{\lambab}{\bar\lambda}

\DeclareMathOperator{\cl}{cl}
\DeclareMathOperator{\deltacl}{\deltabar-acl}
\DeclareMathOperator{\deltadim}{\deltabar-\!\dim}

\DeclareMathOperator{\dacl}{\acl^{\deltabar}}

\newcommand{\eps}{\varepsilon}
\DeclareMathOperator{\rk}{rk}
\DeclareMathOperator{\Jet}{Jet}
\newcommand*{\Jetd}[1]{\Jet_{\delta}^{#1}}

\newcommand{\Ldelta}{L^\delta}
\newcommand{\Ldeltabar}{L^{\bar \delta}}
\newcommand{\monster}{\mathfrak C}
\newcommand{\monstereq}{\mathfrak C^{eq}}

\newcommand*{\intro}[1]{\textbf{#1}}
\newcommand*{\Pa}[1]{\bigl( #1 \bigr)}
\newcommand*{\set}[1]{\{#1\}}
\newcommand*{\abs}[1]{\lvert#1\rvert}

\newcommand*{\card}[1]{\lvert#1\rvert}
\newcommand{\N}{\mathbb{N}}

\newcommand{\R}{\mathbb{R}}

\DeclareMathOperator{\dcl}{dcl}
\DeclareMathOperator{\acl}{acl}
\DeclareMathOperator{\dcleq}{dcl^{eq}}
\DeclareMathOperator{\acleq}{acl^{eq}}
\newcommand*{\tuple}[1]{\langle #1 \rangle}

\newcommand{\av}{\bar a}
\newcommand{\bv}{\bar b}
\newcommand{\cv}{\bar c}
\newcommand{\dv}{\bar d}

\newcommand{\x}{\bar x}
\newcommand{\y}{\bar y}

\newcommand{\deltabar}{\bar \delta}

\DeclareMathOperator{\tp}{tp}

\def\Ind#1#2{#1\setbox0=\hbox{$#1x$}\kern\wd0\hbox to
  0pt{\hss$#1\mid$\hss}\lower.9\ht0\hbox to 0pt{\hss$#1\smile$\hss}\kern\wd0}
\newcommand{\ind}[1][]{\mathop{\mathpalette\Ind{}^{\!\!\!\!\rlap{$\scriptscriptstyle\textnormal{#1}$}\,\,\,\,}}}
\newcommand{\notind}[1][]{\mathrel{\not\mkern-7mu{\ind[#1]}}}
\def\tind{\ind[\th]}
\newcommand{\indd}{\ind[$\deltabar$]}
\newcommand{\indacl}{\ind[M]}
\newcommand{\indacld}{\makebox[1.7\width][l]{$\ind[$M, \deltabar$]$}}
\newcommand{\indf}{\ind[f]}
\newcommand{\indfd}{\ind[f,$\delta$]}
\newcommand{\indtilde}{\tilde\ind}
\DeclareMathOperator{\Uth}{U^{\text{\th}}}

\newcommand{\et}{\ \wedge\ }

\newcommand{\elem}{\equiv}

\def\hyph{\nobreakdash-\hspace{0pt}\relax}

\newcommand{\Wlog}{W.l.o.g\mbox{.}\xspace}
\newcommand{\wloG}{w.l.o.g\mbox{.}\xspace}
\newcommand{\eg}{e.g\mbox{.}\xspace}
\newcommand{\ie}{i.e\mbox{.}\xspace}

\newcommand{\wrt}{w.r.t\mbox{.}\xspace}

\newcommand{\cf}{cf\mbox{.}\xspace}

\newcommand{\Tfae}{T.f.a.e\mbox{.}\xspace}

\newtheorem{lemma}{Lemma}[section]
\newtheorem{thm}[lemma]{Theorem}
\newtheorem{corollary}[lemma]{Corollary}
\newtheorem{mainconjecture}[lemma]{Conjecture}
\newtheorem{proposition}[lemma]{Proposition}
\newtheorem{fact}[lemma]{Fact}
\newtheorem*{fact*}{Fact}
\newtheorem{def-lemma}[lemma]{Definition-Lemma}

\theoremstyle{remark}

\newtheorem{remark}[lemma]{Remark}
\newtheorem{final remark}[lemma]{Final remark}

\newtheorem{claim}{Claim}
\newtheorem*{claim*}{Claim}

\newtheorem{questions}[lemma]{Questions}

\newtheorem{open problem}[lemma]{Open problem}
\newtheorem{conjecture}[lemma]{Conjecture}

\newtheorem{examples}[lemma]{Examples}

\theoremstyle{definition}
\newtheorem{definition}[lemma]{Definition}

\newtheorem{assumptions}[lemma]{Assumptions}

\newenvironment{sentence}[1][]{%
  \begin{list}{}{%
    \setlength\topsep{1ex}%
    \setlength\leftmargin{\parindent}%
  }%
  \item[#1]
 }
 {\end{list}}

\crefname{section}{\S}{\S\S}
\crefformat{section}{\S#2#1#3}
\crefrangeformat{section}{\S\S#3#1#4 to~#5#2#6}
\crefname{subsection}{\S}{\S\S}
\crefformat{subsection}{\S#2#1#3}
\crefrangeformat{subsection}{\S\S#3#1#4 to~#5#2#6}

\crefname{thm}{Thm.}{Thms.}
\crefname{definition}{Def.}{Defs.}
\crefname{proposition}{Prop.}{Props.}
\crefname{corollary}{Cor.}{Cors.}
\crefname{fact}{Fact}{Facts}
\crefname{claim}{Claim}{Claims}
\crefname{assumptions}{Assumptions}{Assumptions}
\crefname{conjecture}{Conjecture}{Conjectures}

\begin{document}

\begin{abstract}
Let T be an algebraically bounded theory. 
We consider the $L(\bar\delta)$-expansions of T by a tuple $\bar \delta$ of derivations (which may be commuting or not).
We investigate the model completion  of either of the above theories, whose
existence has been established in \cite{FT:24}, with particular attention to its
model-theoretic properties, including  $\omega$-stability, simplicity, open core, and elimination of imaginaries.
\end{abstract}



\maketitle


\makeatletter
\renewcommand\@makefnmark%
   {\ \normalfont(\@textsuperscript{\normalfont\@thefnmark})}
\renewcommand\@makefntext[1]%
   {\noindent\makebox[1.8em][r]{\@makefnmark\ }#1}
\makeatother

\tableofcontents

\section{Introduction}

We study the model-theoretic properties of the model completion of the theory of algebraically bounded fields expanded with derivations, continuing the foundational work in \cite{FT:24}.
Algebraically bounded fields, encompassing structures like algebraically closed,
real closed, and p-adically closed fields, provide a rich framework for
exploring the interplay between field-theoretic and model-theoretic algebraic notions.   
Throughout this article,
we denote by  $\K$ a structure that expands a field of characteristic 0.  Recall that $\K$ is  \textbf{algebraically bounded} if the model-theoretic algebraic
closure and the field-theoretic algebraic closure coincide in every structure elementarily equivalent to~$\K$; for more details, examples, and main properties
see~\cite{Dries:89}.  

Let $\K$ be algebraically bounded, 
 $L$ be the language of $\K$, and $T$ be its theory. 
In order to study derivations on $\K$, we denote by  $\deltab = \{\delta_1, \ldots, \delta_k\}$ new unary functions symbols.

Assume that $T$ is model complete.
In \cite{FT:24} we proved that, for every $k \in \N$, the $\Ldb$-expansions of $T$, $\Tdn$ (saying that $\deltab$ is a $k$-tuple of derivations)
and $\Tdb$ (saying that $\deltab$ is a $k$-tuple of \emph{commuting} derivations) have a model completion.
We denote by $\Tdsg$ the model completion of either of the above theories.\\
In the preceding work \cite{FT:24}, we began investigating the model-theoretic properties of $\Tdsg$: we
showed that if $T$ is  stable (resp., dependent), then
$\Tdsg$ is also stable (resp., dependent).

The above results include and extend a long series of results in the literature
(see \cite{FT:24} for a short survey: here we mention \cite{Tressl:05} that shows that large model complete fields admit generic derivations). There are also examples of theories which do not admit generic derivations (see \cite{FT:dexp}).

In this paper, we extend the investigation of $\Tdsg$, 
with a specific focus on model-theoretic properties such as stability, $\omega$-stability, simplicity, and elimination of imaginaries; we also address other questions raised in \cite{FT:24}. Multiple authors have examined these properties in several particular cases; here we consider general setting of  algebraically bounded fields with several generic derivations.
 For instance, in \cite{MS}, the authors demonstrated that $\DCF_{0, m, nc}$, the model completion of the theory of fields with non-commuting derivations, is stable but not $\omega$-stable, eliminates imaginaries, and is uniformly finite. Building upon these results, \cite{Mohamed, SanchM:24} further generalized some of the properties established in \cite{MS}.

In this article, we pay particular attention to the following questions:
\begin{enumerate}
\item When is $\Tdsg$ simple or totally transcendental?
\item Which are the imaginaries of $\Tdsg$? 
\item If $T$ is equipped with a definable topology, is $T$ the open core 
of~$\Tdsg$?
\item Is $\Tdsg$ Uniformly Finite?
\item What is a ``generic'' object?
\item Does $\Tdsg$ have a dimension function in the sense of \cite{Dries:89}?
\item What is the theory of the restriction of $\Tdsg$ where we forget the
derivation but we keep the field of constants?
\end{enumerate}

We will give a full answer to 1) (\cref{sec:simple,sec:w-stable} and \cite{FKM-26}).
In \cref{sec:simple} we prove that if $T$ is simple, then $\Tdsg$ is simple.
\cite{MS} showed that 
the model completion of the theory of fields with non-commuting operators is
simple,
eliminates imaginaries,  is uniformly
finite, and under some additional assumptions it is stable.  

In \cite{FKM-26} the first author together with E. Kaplan and A. Matthews further extend the results exposed here, proving the following:
\begin{thmA}
\begin{itemize}
\item $T$ is rosy iff it is superrosy of thorn rank 1; in this case, $T$ has Geometric Elimination of Imaginaries;
    \item 
    $T$ is simple iff it is supersimple or SU-rank 1;
    \item $T$ is stable iff it is equal to the theory of pure algebraically closed fields, expanded by some constants (we will simply write ``$T = \ACFz$'').
\end{itemize}
\end{thmA}
Combining the above with the results in this paper, we get that $\Tdsg$ is stable iff T is stable iff $T=\ACFz$. In the case the derivations commute, one can improve stability with $\omega$-stability. Furthermore, in the commuting case when T is stable the theory $\Tdsg$ coincides with the theory of differentially closed fields with several commuting derivations $\DCF_{0,m}$; while in the noncommuting case it yields Moosa-Scanlon's theory $\DCF_{0,m,nc}.$


We have some partial answers to 2). We show that if $T$ has (Geometric) Elimination of Imaginaries and it is 
either stable or supersimple, 
then $\Tdsg$ also has (Geometric) EI: \cref{thm:simple-GEI,thm:EI-simple,cor:EI-stable}: together with Theorem A, we get that if $T$ is simple then $\Tdsg$ has GEI. We also show that if $T$ has a definable topology satisfying some general conditions, then $\Tdsg$ eliminates imaginaries
relative to~$T$: see \cref{cor:OC}; we generalize results in~\cite{KP}.  We conjecture that
a stronger result holds than what we manage to prove (see \cref{conj:EI,conj:simple-EI}).

We give a complete answer to 3) (\cref{sec:open} in particular  \cref{cor:OC}).

For 4) we will give a positive answer in a separate paper. 

For 5) we give some possible answers (\cref{sec:Polish}), where we study the Polish space of all derivations on a
countable model $M \models T$ (under some ``bigness'' condition on $M$) and show that
generic derivations form a dense $\Gdelta$-set among all possible derivations.

For 6) we give a complete answer in \cref{sec:ddim}: $\Tdsg$ has a dimension function iff the derivations commute. We also generalize the result in \cite{ELR} on ``coincidence of dimensions'': see
also \cites{GP:12, BMR} for previous results in particular cases.


For 7), in \cref{sec:constant} we study the field of constant $\cfieldK$ of a model
$\tuple{\K, \deltabar} \models \Tdsg$. We
show that $\tuple{\K, \cfieldK}$ is a lovely pair of geometric structures (in
the sense of \cite{BerensteinV:10}), and we
study the definable subsets of $\cfieldK^{n}$.

\smallskip

In \cref{sec:ext,sec:ind} we study independence relations on model of
$T$ and of~$\Tdsg$.
We introduce the independence relation $\indacld$
and prove the fundamental Extension Theorem~\ref{thm:ext} and Independence
Theorem~\ref{thm:ind}: they will make much easier to prove the results in
\cref{sec:simple,sec:open}.
We also characterize the algebraic closure inside models of
$\Tdsg$.

Let us mention here the recent work \cite{PillayPR:25} where the authors study groups
definable in models of $\Tdbg$: they show that any such group can be embedded in
a $T$-definable group.
\medskip






We conclude the paper with several open questions, conjectures, and
announcements of further work.

\section{Preliminaries and conventions}\label{notazioni}
$L$ is a language and $T$ is a $L$-theory (later we will impose additional
 conditions on them).

When we say that $\monster$ is a ``monster model'', we mean that $\monster$ is
a $\lambda$-saturated and $\lambda$-homogeneous model of $T $ for  a ``sufficently big''
cardinal $\lambda > \aleph_{0} + \card L$
(other authors mean that $\monster$ is a proper class, and it is 
$\lambda$-saturated and $\lambda$-homogeneous for every cardinal~$\lambda$).

By ``small'' subset/tuple of $\monster$ we mean  of cardinality less than~$\lambda$.

\smallskip

$\K$ will be an $L$-structure, and $K$ will be either its domain, or, when $\K$
expands a field, its underlying field.

We say that $\K$ is \intro{algebraically bounded} if it expands a field of
characteristic 0 and,  for every subset $A \subseteq K$ (not necessarily a substructure), the model theoretic closure  of $A$ (computed in a sufficiently saturated elementary extension of $\K$) coincides with the field
theoretic algebraic closure of the field generated by $A$ and $\mathrm{dcl}(\emptyset)$.

\subsection{Generic derivations}\label{subsec:prelim-derivation}

This work explores the model theory of differential fields, specifically focusing on expansions of a base theory ~$T$ by derivations.  

Our assumptions for the whole article are the following:

\begin{itemize}
\item $\K$ is a structure expanding a field of characteristic $0$.  
\item $L$ is the language
of $\K$ and $T = Th(\K)$ is its model complete $L$-theory.
\item $F \coloneqq \dcl(\emptyset) \subseteq \K$.
\item $\K$ is algebraically bounded over~$F$.  
\end{itemize}

By algebraic closure we always mean the \(T\)-algebraic closure; in particular, \(\acl\) denotes the \(T\)-algebraic closure, and algebraic independence is understood with respect to \(T\) (equivalently, over \(F\) in the field-theoretic sense).

Under our assumptions, \(\K\) is geometric: in the monster model \(\monster \succ \K\), the operator \(\acl\) satisfies the exchange property, and thus induces a matroid structure.

Moreover, \(\K\) is equipped with a dimension function \(\dim\), assigning to each definable set (with parameters) a natural number and satisfying the axioms of~\cite{Dries:89}.

We also consider the rank function \(\rk\) associated with the matroid \(\acl\): for sets \(V\) and \(B\), \(\rk(V/B)\) is the cardinality of a basis of \(V\) over~\(B\). In particular, if \(X \subseteq M^n\) is definable with parameters \(\bar b\), then
\[
\dim(X) = \max \{ \rk(\bar a / \bar b) : \bar a \in X \}.
\]


Let $\deltabar := \tuple{\delta_{1}, \dotsc, \delta_{k}}$.
Let $\eta_{1}, \dotsc, \eta_{k}$ be  derivations on~$F$.
We obtain two theories:

\begin{description}
\item[$\Tdb$] the expansion of $T$ saying that the $\delta_{i}$ are derivations which
commute with each other, i.e. for every $i$, $j \leq k$,
$\delta_{i} \circ \delta_{j} = \delta_{j} \circ \delta_{i}$, and that $\delta_{i}$ extends $\eta_{i}$ for $i \leq k$;
\item[$\Tdn$] the expansion of $T$ saying that the
$\delta_{i}$ are derivations without any further conditions and that $\delta_{i}$ extends $\eta_{i}$ for $i \leq k$.
\end{description}
Both theories have a model completion (remember that $T$ is model complete) 
(see \cite{FT:24}).
For convenience, we use $\Tdsg$ to denote either of the model completions,
both for commuting and non-commuting tuples of derivations.\\
$\tuple{\monster, \deltab}$ is a monster model of $\Tdsg$.\\
$\tuple{\K, \deltab}$ is some model of $\Tdsg.$

\subsection{Notations and know results}\label{notation}

We introduce some notation that we use in the sequel.\\
We denote by $\Gamma$ the free commutative (or non commutative) monoid generated by $\deltabar,$ we can define the
canonical partial order~$\preceq$ given by $\beta \preceq \alpha\beta$, for all $\alpha, \beta \in \Gamma$.

\begin{remark}
If $\Gamma$ is the free non commutative monoid generated by $\deltabar$ then $\preceq$ is a well-founded partial order on $\Gamma$, but it is not a well-partial-order
(\ie, there exist infinite anti-chains).


\begin{enumerate}
\item $\emptyset$ (\ie, the empty word, corresponding to the identity function on $\K$)
 is the minimum of $\Gamma$;
\item If $\alpha \preceq \beta$, then $\gamma \alpha \leq \gamma \beta$ and $\alpha \gamma \leq \beta \gamma$.
\end{enumerate}
\end{remark}

\begin{remark} If $\Gamma$ is the free commutative monoid generated by $\deltabar$, with the
canonical partial order~$\preceq$,  it is isomorphic to $\N^k$.
\end{remark}


Fix $n \in \N$, $\x = \tuple{x_{1}, \dotsc, x_{n}}$, 
and denote by $\Gamma_{n} \coloneqq
\set{\gamma x_{i}: \gamma \in \Gamma, i \leq n}$: we consider $\Gamma_{n}$ 
as a set of some $\Ldb$-terms, as a set of functions from $K^{n}$ to~$K$, and as a set
of indices. 
For example $\delta_1\delta_2x_3$ is in $\Gamma_{3}$ but   $\delta_1\delta_2x_3 + \delta_1x_1$ is not. 
Given $\bar a \in \K^n$, we denote by $\Jet(\av) := \tuple{\gamma \av : \gamma \in \Gamma_n}$. 

Moreover, we denote by $\x_{\Gamma} $ a set of formal variables indexed by~$\Gamma_{n}$: that
is, for every $f \in \Gamma_{n}$ we have a variable $\x_{f}$.
Given an $L(A)$-formula $\alpha(\x_{\Gamma})$
we denote by $\alpha(\Gamma)$ the $\Ldb(A)$-formula where we replaces each
occurrence of a variable $\x_{f}$ with the corresponding term~$f(\x) \in \Gamma_{n}$.



Let $J \subseteq \Gamma_{n}$ and $A \subset K$ be a ``small'' subset of $K$ with $\deltabar A \subseteq A$.
Let $p(\x)$ be a complete $\Ldb$-type over $A$ (in $n$ variables).
Define
\[
p_{J} \coloneqq \set{ \alpha: \alpha(\x_\Gamma)\ L(A)\text{-formula s.t. } \alpha(\Gamma) \in p} \in
S_{L}^{J}(A).
\]
We say that some type $q \in S^{J}_{L}(A)$ is $\deltabar$-compatible if there exists
 $s \in S^{n}_{\Ldb}(A)$ such that $q = s_{J}$.

\medskip


Many of the model-theoretic properties of $T$ are inherited by~$\Tdsg$. The following results from  \cite{FT:24}, which plays a crucial role in our subsequent analysis, is particularly useful:

\begin{thm}
For  every $\av$ tuple in $\monster$ and $B$
subset of $\monster$, the $\Ldb$-type of $\av$ over $B$ is uniquely determined by the $L$-tuple
of $\Jet(\av)$ over $\Jet(B)$.
\end{thm}

\begin{thm}\label{thm:strong-EQ}
\begin{enumerate}
\item If $T$ eliminates quantifiers, then $\Tdsg$ eliminates quantifiers.
\item For every $\Ldeltabar$-formula $\alpha(\x)$ there exists an $L$-formula $\beta(\x_\Gamma)$
such that
\[
\Tdsg \models \forall \x\ \Pa{\alpha(\x) \leftrightarrow \beta(\Gamma)}.
\]
\end{enumerate}

\end{thm}

\section{Elimination of  imaginaries}

Elimination of imaginaries is a powerful concept in model theory that simplifies
arguments by allowing us to work directly with elements of the model, rather
than with equivalence classes of definable sets. 


We remind  what Elimination of Imaginaries and some of its variants are.
\begin{definition}
$T$ eliminates imaginaries ($T$ has EI) if every imaginary is interdefinable with a real tuple: i.e., for every
imaginary $e$ there is a real tuple $a$ such that $\dcleq(e) = \dcleq(a).$
\end{definition}

\begin{definition}
$T$ weakly eliminates imaginaries ($T$ has WEI), if for every imaginary $e$ there is a real tuple $a$ such that $e$ is
definable over $a$ and a is algebraic over~$e$: 
that is, $e \in \dcleq(a)$ and $a \in \acleq(e).$
\end{definition}

\begin{remark}
Equivalently, a theory $T$ admits weak elimination of imaginaries (WEI) in
the sense of B. Poizat (see \cite{CasanF:04,Yoneda:22}) if, for any $\phi(\overline x, \overline a) \in  L(\overline a)$,  we have the
smallest algebraically closed set $B$ such that  $\phi(\overline x, \overline a) $
is definable over~$B$.
\end{remark}

\begin{definition}
$T$ has geometric eliminates imaginaries ($T$ has \GEI), if for every imaginary $e$ there is a real tuple $a$ such that $e$ is
algebraic over $a$ and $a$ is algebraic over $e$: that is, $e \in \acleq(a)$ and 
$a \in \acleq(e)$.
\end{definition}

\begin{remark}
By definition, it is clear that EI implies WEI, and that WEI implies GEI. The converse implication, namely that WEI implies EI, holds when working over a field as is the case here (see \cite{MMP}).
\end{remark}

\begin{remark} \label{WEIimpliesEI}In any structure containing at least two constants, as in our
case, the definition of EI is equivalent to the Uniform Elimination of
Imaginaries (see \cite{CasanF:04} for definition and proof).
\end{remark}

\bigskip

Regarding generic derivations, we formulated the following:
\begin{conjecture}[\cite{FT:24}]\label{conj:EI}
Let $T$ be algebraically bounded.
Then, $\Tdsg$ has  elimination of imaginaries modulo $T^{eq}$.

In particular: if $T$ has \GEI, then $\Tdsg$ also has \GEI, and if $T$ has EI, 
then $\Tdsg$ also has EI.
\end{conjecture}

A few particular cases were already known, when $T$ is one of the following:
\begin{itemize}
\item 
$\ACFz$: see \cite{mcgrail,MS};
\item $\RCF$: 
see \cite{FK, KP} for a
proof based on M. Tressl's idea, see also \cites{bkp, Point} for different proofs.
\end{itemize}

We will prove two cases of the above conjecture:
\begin{itemize}
\item  When $T$ is simple: by Theorem~A $T$ has \GEI;
we then employ
a result in \cite{Yoneda:09} to show
that $\Tdsg$ also has \GEI (see \cref{sec:simple-GEI}); if $T$ has EI, we 
adapt a technique in \cite{HC:99} to show that $\Tdsg$ also has EI
(see \cref{sec:simple-EI}). 
\item When $T$ has a suitable definable topology: we  employ a technique by
M.~Tressl to show that $\Tdsg$ has elimination of imaginaries modulo $T^{eq}$
(see \cref{sec:open}).
\end{itemize}

\section{Independence relations}\label{sec:independence-relations}

To establish the primary results of this paper, we remind the concept of independence relations and revisit key theorems related to this notion.  We fix a monster model $\monster \models T$.

Let $\ind$ be a ternary relation on small
subsets of~$\monster$.
We mostly use the definitions and nomenclature from \cite{Adler}: in particular, we are
interested in the cases when $\ind$ is either a strict independence relation, or
the relation $\indacl$ defined in \cite[Def.~1.26]{Adler}, 
or Shelah-forking~$\indf$.

We remind for completeness the definition $\indacl$  introduced in \cite{Adler}.
\begin{definition}\label{def:Mind}
The relation $\indacl$ (M-dividing independence) is defined as:

$A \indacl_C B$ iff for any $C'$ s.t $C \subseteq C'  \subseteq  \acl(BC)$ then $ \acl(AC') \cap \acl(BC') = \acl C'$.
\end{definition}

For the remainder of this section, $\ind$ is an independence relation 
on~$\monster$.

\begin{remark}
If $A \ind_{C} B$, then $A \ind_{C} \acl(BC)$.
\end{remark}
\begin{proof}
By Existence, there exists $D \elem_{BC}{\acl(BC)}$ s.t. $A \ind_{BC} D$: notice
that $D = \acl(BC)$, and therefore $A \ind_{BC} \acl(BC)$.
By transitivity, $A \ind_{C} \acl(BC)$.
\end{proof}

Notice that we don't assume that $T$ has some form of elimination of
imaginaries, nor that $\ind$ extends to $\monster^{eq}$.
However, the setting in \cite{Adler} often assumes implicitly to work in
$T^{eq}$: to avoid confusion, we will  introduce some additional nomenclature,
following in part \cite{Yoneda:09}.

\begin{definition}\label{def:strict}
We say that an independence relation $\ind$ is:
\begin{description}
\item[real-strict] if, for every $a \in \monster$ and $B \subset \monster$,
$a \ind_{B} a$ iff $a \in \acl(B)$;
\item[strict] if $\ind$ extends to an independence relation
on $\monstereq$ (which we also denote by~$\ind$),
and for every $a \in \monstereq$ and $B \subset \monstereq$, $a \ind_{B} a$ iff
$a \in \acleq(B)$;
\item[real-canonical] if it is \realstrict and, for every $A, B, C, D \subset \monster$ such that $B \supseteqq C, D$ 
\[
A \ind_{C} B \et A \ind_{D} B \implies A\ind_{\acl(C) \cap \acl(D)} B;
\]
\item[canonical] if it is \eqstrict and, for every $A, B, C, D \subset \monstereq$, such that  $B \supseteqq C, D$ 
\[ 
A \ind_{C} B \et A \ind_{D} B \implies A\ind_{\acleq(C) \cap \acleq(D)} B
\]
\end{description}
\end{definition}

Notice that if $\ind$ is \eqstrict, then it is \realstrict.
If either $\monster = \monstereq$, or $\monster$ has \GEI (Geometric Elimination of Imaginaries), then $\ind$ is
\realstrict iff it is \eqstrict;
similarly, under the same assumption, $\ind$ is \realcanon iff it
is \eqcanon. 

\begin{fact}
\begin{enumerate}
\item 
There exist rosy theories with EI (even o-minimal) such that $\tind$
(thorn-forking independence)
is not (real)-canonical: \cite[Appendix A.3]{Adler}
explains an example from \cite{LP:93}; see also \cite[Remark~3.5]{Yoneda:09}.
\item
There exist a theory with WEI (the ``integral Urysohn
space'') and more than one strict independence relation: see \cite{Conant-16}.
\end{enumerate}
\end{fact}

\begin{fact}[{\cite[Lemma~3.1]{Yoneda:09}}]\label{fact:Yoneda}
If $\ind$ is \eqstrict and \realcanon, then $\monster$ has \GEI.
\end{fact}

\begin{corollary}\label{cor:Yoneda}
Let $\ind$ be \realstrict and \realcanon.
\Tfae:
\begin{enumerate}
\item $\monster$ has \GEI; 
\item $\ind$ is \eqstrict;
\item $\ind = \tind$.
\end{enumerate}
\end{corollary}

\begin{lemma}\label{lem:strict-M}
Let $\ind$ be a  \realstrict 
independence relation on~$\monster$.
If $A \ind_{C} B$, then $A \indacl_{C} B$.

In particular, 
if $T$ is simple and $A \indf_{C} B$, then $A \indacl_{C} B$.
\end{lemma}
\begin{proof}
Let $A, B, C, C'$ be as in \cref{def:Mind}.
Assume that $A \ind_{C} B$.
Let $d \in \acl(AC') \cap \acl(BC')$.
We need to show that $d \in \acl(C')$.
\begin{multline*}
A \ind_{C} B \implies A \ind_{C} \acl(BC) \implies A \ind_{C} BC'
\implies 
\\
\implies  A \ind_{C'} B \implies \acl(AC') \ind_{C'} \acl(BC')
\implies d \ind_{C'} d.
\end{multline*}
Since $\ind$ is \realstrict, the latter implies that $d \in \acl(C')$.
\end{proof}

\begin{fact}[{\cite{BuechPW:00},\cite[\S3.1]{Adler}}]\label{fact:canonical}
\begin{enumerate}
\item If $\monster$ is either stable or super-simple, then it has
 Elimination of HyperImaginaries (EHI);
\item If $\monster$ is simple with EHI, then $\indf$
is \eqcanon;
\item If $\ind$ is \eqcanon, then $\ind = \tind$.
\end{enumerate}
\end{fact}
For more on EHI in (simple) theories, see e.g. \cite{Casanovas:11, PalacW:13, BuechPW:00}; 
there are no known examples of simple theories that do \emph{not}
eliminate hyperimaginaries; for us the relevant implication is the following fact:
\begin{fact}\label{cor:canonical}
If $\monster$ is either stable or super-simple, then $\indf$ is
\eqcanon and $\indf = \tind$.
\end{fact}

\begin{definition}
$\ind $satisfies ``Independence over Models'' if:
\begin{sentence}
For every $\av, \bv \in \monster^{n}$, $\av', \bv' \in \monster^{m}$,
for every $M \prec \monster$ with $M$ small, if
\[
\av \ind_{M} \bv, \qquad \av' \ind_{M} \av, \qquad \bv' \ind_{M} \bv, \qquad \av'
\elem^{L}_{M} \bv'
\]
then there exists $\cv \in \monster^{m}$ such that
\[
\cv \ind_{M} \av \bv, \qquad \cv \elem^{L}_{M \av} \av', \qquad \cv \elem^{L}_{M
  \bv} \bv'.
\]
\end{sentence}
$\ind $ satisfies ``Independence over Algebraically
Closed (real) Substructures'' (IACS) if the above is true when we let $M$ vary
among algebraically closed substructure of~$\monster$.
We say that $T$ satisfies \IACS if $\indf$ satisfies it.
\end{definition}

\begin{fact}[Kim-Pillay: see {\cite{KimP:98}, \cite[Thm.~29.11, Thm.~29.13]{TZ}}]\label{fact:KP}
Let $\ind$ be an independence relation on 
$\monster$ that
satisfies Independence over Models.
Then, $\monster$ is simple and $\ind$ is Shelah-forking
(on~$\monster$).

Conversely, if $\monster$ is simple, then $\indf$ is a strict
independence relation satisfying Independence over Models.
\end{fact}

\begin{fact}[\cite{BuechPW:00}]\label{fact:IACS}
If $\monster$ is simple with EHI and EI, then $\indf$ satisfies \IACS.
\end{fact}

$\ind$ can be extended to $L$-types: given $A \subseteq B$ small subsets of
$\monster$ and an $L(B)$-type (in a small number of variables) $q(\y)$, we write
\[
q \ind_{A} B
\]
if for \emph{some} $\cv \in \monster^{\card \y}$ realizing   $q(\y)$ we have that
\[
\cv \ind_{A} B.
\]
Notice that the above is equivalent to:

``For \emph{every} $\cv \in \monster^{\card \y}$ realizing  $q(\y)$ we have that
$\cv \ind_{A} B$''.

\medskip


\subsection{Independence relations on algebraically bounded structures}

We analyze now the case when $L$ extends the language of rings and $T$
is algebraically bounded.

\begin{lemma}
Let $A, B, C \subseteq \monster$.
\Tfae:
\begin{enumerate}
\item $A \indacl_{C} B$;
\item for every $A' \subseteq A$, if $A'$ is algebraically independent over $CF$ (where $F = dcl(\emptyset)$), then it
is still algebraically independent over $BC$;
\item for every $A' \subseteq A$, if $A'$ generates $A$ (in the sense of the
matroid $\acl$) over $BC$, then it
generates $A$ over $B$;
\item\label{en:ind-sum} for every $A' \subseteq A$ and $B' \subseteq B$,
$\td(A'B'/C) = \td(A'/C) + \td(B'/C)$.
\end{enumerate}
Moreover, $\indacl$ is a \realstrict independence relation on $\monster$ satisfying
the Strong Finite Character condition.
\end{lemma}
\begin{proof}
For the ``moreover'' part, notice that \ref{en:ind-sum} implies that
$\indacl$ is symmetric, and the conclusion follows by 
\cite[\Thms 2.39 and Exercise~2.6]{Adler}.
\end{proof}

\begin{fact}\label{fact:GEI-rosy}
If $\monster$ is algebraically bounded and has \GEI, then $\indacl = \tind$ and
$\monster$ is rosy, i.e.  $\tind$  is an independence relation on
$\monster^{eq}$; moreover, it is super-rosy of \th-rank 1.%
\footnote{Later, we will prove a kind of converse to this fact.}
\end{fact}
\begin{proof}
Since $\monster$ has \GEI, $\indacl$ is a strict independence relation on
$\monster^{eq}$, and therefore $\monster$ is rosy.
For more details see \cite{Adler}: in particular, \cite[Def~1.28]{Adler}.

$\monster$ is then super-rosy of rank 1 because $\Uth(\av/B) =\td(\av/ FB)$,
where $\Uth$ is the \th-rank,  $\td$ is the transcendence degree, and
$F$ is the model-theoretic algebraic closure of the empty set.
\end{proof}

\begin{remark}\label{rem:simple-GEI-ind}
Assume that $\monster$ is algebraically bounded, with \GEI, simple, and $\indf$ is
canonical.
Then, $\indacl = \indf$.
\end{remark}
\begin{proof}
Since $\indf$ is canonical, $\indf= \tind$, and by the Fact \ref{fact:GEI-rosy} we have the conclusion. 
\end{proof}

\section{The Extension Theorem}
\label{sec:ext}

The results presented in this section are of independent interest and will be utilized in the sequel. 
For now, $K$ is some field of characteristic~$0$.
\begin{definition}
We denote by $\indacl$ the ternary relation on subsets of $K$ induced by $\acl$ over $F$:
\[
A \indacl_{C} B
\]
iff, for every $\av$ finite subset of $A$, if $\av$ is algebraically independent
over~$CF$, then $\av$ remains algebraically independent over $CFB$.%
\footnote{When $\K$ is algebraically bounded, the above definition coincides
  with the one in \cref{sec:independence-relations}.}
\end{definition}

The following lemma is well-known; we provide a proof for completeness.

\begin{lemma}[Amalgamation] \label{amalgamation}
Let $B_0$, $B_1$, $B_2$ be subrings of $K$ containing~$F$.
Let $\deltabar_i$ be $k$-tuples of derivations on $B_i$, $i=0, 1, 2$.
Assume:
\begin{enumerate}
\item $\tuple{B_0, \deltabar_0}$ is a common substructure of 
$\tuple{B_1, \deltabar_1}$ and $\tuple{B_2, \deltabar_2}$;
\item \[B_{1} \indacl_{B_{0}} B_{2}.\]
\end{enumerate}
Then, there exists a $k$-tuple of derivations $\deltabar$ on $K$ which extends all the~$\deltabar_i$.
Moreover, if each  $\deltabar_i$ commutes, then $\deltabar$ commutes.
\end{lemma}
\begin{proof}
Let $\bv_i$ be a transcendence basis of $B_i$ over $B_0$, $i = 1,2$.
By 2), $\bv_1$ and $\bv_2$ are disjoint, 
and $\bv_1 \cup  \bv_2$  is algebraically independent over $B_0$.
Thus, there exists a  $k$-tuple of derivations  $\deltabar$ on $K$ 
extending $\deltabar_{0}$  and such that
\[
\deltabar(\bv_i) = \deltabar_i(\bv_i), \qquad i= 1, 2.
\]
But then $\deltabar$ extends also $\deltabar_1$ and $\deltabar_2$.
If moreover each $\deltabar_i$ commutes, let $\cv$ be a transcendence basis of
$K$ over 
$B_0 \cup B_1 \cup B_2$.
Let $\deltabar$ be the unique $k$-tuple of derivations on $K$  such that:
\begin{enumerate}
\item $\deltabar$ extends $\deltabar_{0}$; 
\item $\deltabar(\bv_i) = \deltabar_i(\bv_i)$, $i= 1, 2$; 
\item $\deltabar (\cv )$ = 0.
\end{enumerate}
Then, $\deltabar$ commutes.
\end{proof}

Now $\tuple{\K, \deltabar}$ is a monster model of $\Tdsg$.

The following  technical result has several applications: we will see some
later;  in \cite[Lemma~5.8]{PPP:23} the authors prove a version of it in order to study definable groups.

\begin{thm}[Extension]\label{thm:ext}
Let $A \subseteq B$ be ``small'' subsets of $K$ with $\deltabar A \subseteq A$ and 
$\deltabar B \subseteq B$.
Let $p(\x)$ be a complete $\Ldb$-type over $A$ (in $n$ variables).
Let $q \in S_{L}^{J}(B)$ be some extension of $p_{J}$.
If
\[
q \indacl_{A} B,
\]
then $q$ is $\deltabar$-compatible (see Subsection \ref{notation}). 
\end{thm}
\begin{proof}
\Wlog, $A$ and $B$ are also subrings of~$K$.
Let $\av \in K^{n}$ be a realization of $p$ and $\bv \in K^{J}$ be a realization
of~$q$.
Since $\Jet(\av)$ and $\bv$ have the same $L$-type $q$ over $A$, there exists
$\phi$ an $L$-automorphism of $\K$ fixing $A$ point-wise and mapping $\bv$ to
$\Jet(\av)$.

Define $\epsb \coloneqq \phi^{-1} \circ \deltabar \circ \phi$.
Let
\begin{itemize}
\item $\deltab_{0}$ be the restriction of $\deltab$ to $A$, 
\item $\deltab_1$ be the restriction of $\deltab$ to $B$, 
\item  $\deltab_{2}$ be the restriction of $\epsb$ to $A[\bv]$.  
\end{itemize}
Notice that $\epsb$ is a tuple of derivations on $K$ which extends
$\deltabar_{0}$ and such that $\epsb$
 commutes if $\deltabar$ commutes.
Thus, by the Amalgamation Lemma, there exists a $k$-tuple 
$\lambab$ 
of derivations on $K$
which extends $\deltab_0$, $\deltab_{1}$, $\deltab_{2}$ and such that $\lambab$ commutes if $\deltab$ commutes.
Thus, $\tuple{\K, \lambab} \models \Tds$ ($\lambab$ is not generic, in general).
Let $ \bv_0 \in K^{n}$, interpreted in  $\tuple{\K, \lambab} $, we have that $\Jet( \bv_0) = \bv$, and therefore 
$\Jet( \bv_0)$ satisfies~$q$.
Thus, the partial $\Ldb(B)$-type
\[
u(\x) \coloneqq \Tdsg \cup  \Diag_{\Ldb}(B) \cup q(\Jet(\x))
\]
is consistent.
Since $\Kdb$
 is existentially closed and $\card B^{+}$-saturated, 
there exists $\cv \in K^{n}$ satisfying~$u$. 
\end{proof}

\subsection{Relation with the framework of \texorpdfstring{\cite{SanchM:24}}{LM25}}
It is worth noting that, in view of Lemma~\ref{amalgamation}, the theory 
$T^{\bar\delta}$ with respect to $(T, \indacl)$ fits into the framework 
of derivation-like expansions introduced in \cite{SanchM:24}. 

Indeed, Lemma~\ref{amalgamation} yields the amalgamation property required 
in part~(1) of the definition of derivation-like expansions. Moreover, 
since $\K$ is algebraically bounded, algebraic closure defines an 
independence relation, and, in characteristic $0$, derivations extend 
uniquely to algebraic extensions. This immediately implies part~(2) of 
the definition.

Furthermore, when $T = \mathrm{Th}(\K)$ is simple, forking independence 
implies algebraic independence, and hence $T^{\bar\delta}$ is derivation-like 
also with respect to $(T, \indf)$.

As a consequence, several of our results can alternatively be obtained from the general theory developed in \cite{SanchM:24}, 
and our Theorems \ref{thm:ext}, \ref{thm:ind}, \ref{thm:simple}, \ref{thm:simple-d} (as well as Corollary \ref{cor:acld}) can be recovered from their corresponding general results in that setting.

Nevertheless, we have chosen to include self-contained proofs of these results. This is partly for the convenience of the reader, and partly to maintain a uniform approach and notation throughout the paper, avoiding the need to repeatedly translate between our setting and the abstract framework of \cite{SanchM:24}.

\section{Algebraic closure and independence relations}
\label{sec:ind}
In this section we examine the relations between algebraic closure, independence
relations, and (geometric) elimination of imaginaries.
The results in this section will be used in
\cref{sec:simple,sec:ddim}.

We fix $\tuple{\monster; \deltabar}$ monster model of $\Tdsg$.

Let $\ind$ be some ternary relation on subsets of $\monster$. 
We define the following ternary relation on subsets of $\monster$.
\[
A \indd_{C} B \iff \Jet(A) \ind_{\Jet(C)} {\Jet(B)}.
\]
Notice that we are not assuming that $T$ eliminates imaginaries,
or that $\ind$ extends to $\monster^{eq}$).

\begin{thm}\label{thm:ind}
Assume that $\ind$ is a \realstrict independence relation on~$\monster$
(see \cref{def:strict}).
Then,
$\indd$ is a \realstrict independence relation on $\tuple{\monster; \deltabar}$.


If moreover $\ind$ is \realcanon, then $\indd$ is \realcanon.

\end{thm}
\begin{proof}
In \cite{Adler}, a list of axioms for independence relations is presented. The
author subsequently demonstrates that the Extension Axiom is equivalent to the
conjunction of the Existence Axiom and the Symmetry Axiom. Notably, among the
axioms for independence ($\indd$) outlined in \cite{Adler}, only the Existence
Axiom requires non-trivial verification. 
Thus, let $\av, \bv, \cv$ be small tuples in $\monster$.
Since $\ind$ is an independence relation, there exists a small tuple
$\dv \subset \monster$ s.t. $\dv \ind_{\Jet (\cv)} \Jet( \bv)$.
Let $p \coloneqq \tp^{\Ldb}(\Jet (\av)/ \Jet (\cv))$ and $q \coloneqq \tp^{L}(\dv /
\Jet(\cv) \Jet (\bv))$.
By assumption, $q$ is a complete $L$-type extending $p$, and 
\[q \ind_{\Jet (\cv)} \Jet(\bv).\]
Since $\ind$ is \realstrict, \cref{lem:strict-M} implies that \[q \indacl_{\Jet (\cv)} \Jet(\bv)\] and therefore,
 by the Extension Theorem, there exists $\av' \subset \monster$ s.t.\
$\Jet (\av')$ realizes~$q$: hence,
\[
\av' \elem^{\Ldb}_{\cv} \av \et \av' \indd_{\cv} \bv,
\] 
proving that $\indd$ satisfies Existence.



\smallskip

The fact that $\indd$ is \realstrict is clear.

\smallskip
Finally, assume that $\ind$ is \realcanon.
Assume that 	 
\[
\cv  \indd_{E\bv}  \av \et \cv  \indd_{E\av}  \bv.
\]
We have to prove that
\[
\cv  \indd_{\dacl(E \av) \cap \dacl(E \bv)}  \av \bv.
\]
\Wlog, we may assume that $E = \dacl E$.
Thus, we have 
\[
\Jet(\cv)  \ind_{E \Jet(\bv)}  \Jet(\av) \et \Jet(\cv)  \ind_{E\Jet(\av)}  \Jet(\bv).
\]
Since $\ind$ is \realcanon, we have
\[
\Jet(\cv)  \ind_{\acl(E \Jet(\bv)) \cap \acl(E \Jet(\bv))} \Jet(\bv).
\]
$\acl(E \Jet(\bv)) = \dacl(E  \bv)$, and similarly for $\av$,
we are done.
\end{proof}


In particular, we can choose $\ind := \indacl$,  and obtain that the induced relation $\indacld$ 
is a \realstrict independence
relation on $\tuple{\monster; \deltabar}$.

\begin{corollary}[{\cf \cite[Cor.1.4]{PillayPR:25}, \cite[Cor.3.16]{SanchM:24}}]
\label{cor:acld}
1) The algebraic closure on $\tuple{\monster; \deltabar}$ is given by
\[
\acl(\Jet(A)),
\]
for every $A \subseteq \monster$.%
\footnote{$\Tdsg$ is differentially algebraically bounded: see 
\cite[\Defs 4.1, \Rems 4.21]{Wang:25}; \cf \cite[\S2.25]{Dries:89}.}

2) If $\tuple{\monster; \deltabar}$ has \GEI, 
then  it is rosy, and $\indacld$ coincides with \th-forking $\tind$ (for $\tuple{\monster; \deltabar}$).
\end{corollary}
\begin{proof}
Let $\dacl$ be the algebraic closure according to $\Tdsg$.

1) If $a \in \dacl B$, then, since $\indacld$ is an independence relation,
\[
a \indacld_{B} a,
\]
and in particular
\[
a \indacl_{\Jet(B)} a
\]
that is $a \in \acl(\Jet(B))$.

Conversely, it is clear that $\acl(\Jet(B))$ must be contained in $\dacl(B)$.

2) Since $\indacld$ is a \realstrict independence relation on $\tuple{\monster, \deltab}$, we have that the latter is real-rosy; if
moreover $\tuple{\monster, \deltab}$ has \GEI, it is rosy (see \cite{Adler} for definitions
and proofs: in particular, \cite[Def~1.26]{Adler}).
\end{proof}


\subsection{GEI in algebraically bounded
structures}

\begin{proposition}\label{prop:RIP-acl}
Let $K$ be any field of characteristic $0$.
Define $A \ind_{C} B$ if, for every $X \subseteq A$ which is algebraically independent
over $C$ (in the field-theoretic sense), $X$ remains algebraically independent
over $BC$.
Then, $\ind$ is \realcanon.
\end{proposition}
\begin{proof}
Let $\Kt$ be the field-theoretic algebraic closure of $K$.
Given $C \subseteq K$, let $\tilde C$ be the field-theoretic algebraic closure of
$C$ inside $\Kt$, and $\acl(C)$ be the field-theoretic algebraic closure
of $C$ inside $K$.

Let $D \coloneqq C_{1} \cap C_{2}$.

\begin{claim}\label{cl:absolute-acl-intersection}
\[
\tilde{C_{1}} \cap \tilde{C_{2}} = \widetilde {\acl(C_{1}) \cap  \acl(C_{2})}.
\]
\end{claim}
First, \wloG we may assume that $C_{i} = \acl(C_{i})$, for $i = 1, 2$.
It is clear  that
\[
\tilde D \subseteq \tilde{C_{1}} \cap \tilde{C_{2}};
\]
we have to prove the opposite inclusion.
Let $C$ be either $C_{1}$ or $C_{2}$.
Let $a \in \tilde C_{1} \cap \tilde C_{2}$ and $a_{1}, \dotsc, a_{\ell}$ be the  conjugates of $a$ over  $K$.
Then, $a_{j} \in \tilde C$, $j= 1, \dotsc, \ell$.
Let $\sigma_{t}(x_{1}, \dotsc, x_{\ell})$ be one of the symmetric polynomials: we have
that
$\sigma_{t}(a_{1}, \dotsc, a_{\ell}) \in \tilde C \cap K = C$.
Thus, if $p(x) \in K[x]$ is the minimal polynomial of $a$ over $K$, we have
$p \in C[x]$, and therefore $p \in D[x]$.

Thus, $a \in \tilde D$, proving the claim. 

\smallskip

Let $\indtilde$ be the analogue of $\ind$ inside $\Kt$.

\begin{claim}
$\indtilde$ is  \realcanon. 
\end{claim}
In fact,
$\indtilde$ coincides with $\indf$ on $\Kt$ seen as a pure field, and, 
since $\Kt$ (using the fact that it is an algebraically closed field) is stable
with EI, $\indf$ is \realcanon: see \cite[\S3.1]{Adler}.
\smallskip

Assume now that $C_{i} \subseteq B$ for $i = 1, 2$, and
$A \ind_{C_{i}} B$ and $C_{i} \subseteq B$, for $i = 1, 2$.

We have $A \indtilde_{C_{i}} B$ for $i = 1, 2$.
Since $\indtilde$ is \realcanon, we have that
\[
A \indtilde_{\tilde C_{1} \cap \tilde C_{2}} B,
\]
and therefore, by \cref{cl:absolute-acl-intersection},
\[
A \indtilde_{\tilde D} B,
\]
which is equivalent to $A \ind_{D} B$.
\end{proof}

\begin{corollary}\label{cor:indacl-RC}
$\indacl$ is \realstrict and \realcanon.\\
$\indacld$ is also \realstrict and \realcanon.
\end{corollary}

\begin{proof} This is an immediate consequence of \cref{thm:ind} and 
\cref{prop:RIP-acl}.

\end{proof}

Thus, by \cref{cor:Yoneda}, we have:
\begin{thm}\label{thm:GEI-rosy-eq}
\Tfae:
\begin{enumerate}
\item $\monster$ has \GEI; 
\item $\indacl$ is \eqstrict;
\item $\indacl = \tind$;
\item $\monster$ is super-rosy of \th-rank 1.
\end{enumerate}
\end{thm}

\begin{thm}\label{thm:delta-GEI-rosy-eq}
\Tfae:
\begin{enumerate}
\item $\tuple{\monster; \deltabar}$ has \GEI; 
\item $\indacld$ is \eqstrict;
\item $\indacld$ is equal to \th-forking on $\tuple{\monster; \deltabar}$.
\end{enumerate}
\end{thm}

\begin{conjecture}\label{conj:simple-EI}
Assume that $\monster$ is has  \GEI.
Then, $\tuple{\monster; \deltabar}$ also has \GEI.
\end{conjecture}



\begin{proposition}\label{prop:HC-claim}
Let $\monster$ be a monster model of some theory. Assume that $\ind$ is an independence
relation on $\monster$ satisfying the following conditions:
\begin{enumerate}[Romanenum]
\item $\ind$ is \realcanon;
\item\label{en:controlled-rank}
$\ind$ is ``controlled by a rank'': 
 there exist a linearly ordered set $\tuple{I ; \leq}$ and a function $\rho$
associating an element of $I$
to each type (over some small subset of the monster model $\monster$ in finitely
many real variables), 
such that:
\begin{enumerate}
\item 	$\rho$ is invariant under automorphisms,
\item	$\rho$ for every $\av$, $B$, $C$,
	$\rho(\av/B) \leq \rho(a/BC)$, 
	with equality iff  $\av   \indf_{B}  C$,
\item\label{en:USC} for every non-empty  set $X$ which is type-definable
over~$A$, $\rho(\x/A)$ attains a maximum on~$X$ (e.g., $\rho$ is upper 
semi\hyph continuous).
\end{enumerate}
\end{enumerate}
Let $e \in \monstereq$, where $e = f(\av)$ for some (real finite tuple)
$\av \in \monster^{\ell}$ and some function $f$ which is definable without parameters.
Let $E \coloneq \acleq(e) \cap \monster$ and $P$ be the set of realization of $\tp(a/E)$ (inside~$\monster$). Then there exists $\cv \in P$ such that $f (\cv) = e$ and $\av \ind_{E} \cv$.
\end{proposition}

\begin{proof}
\ref{prop:HC-claim} follows the proof in \cite[(1.10)]{HC:99} and \cite[\Thms 5.12]{MS}.
As noticed in \cite{HC:99}, Neumann's Lemma implies that there exists
$\bv$ conjugate of $\av$ over $E \cup \set e$ such that
\begin{equation}\label{eq:1}
\acl(E\av) \cap \acl(E\bv) = \acleq(E e) \cap \monster  = E
\end{equation}
(see \cite{EH:93}).
Let $Q$ be the set of $\bv \in P$ satisfying \eqref{eq:1}.
\begin{claim}
$Q$ is type-definable over $E \av$.
\footnote{This step was not clear to us: thanks to Itay Kaplan for explaining it.}
\end{claim}
Since $P$ is type-definable, it suffices to show that the set
\[
X \coloneqq \set{\x \in \monster^{\ell}: \acl(E\av) \cap \acl(E\x) \supsetneq E}
\]
is \ordefinable.
For every $d \in \acl(E\av) \setminus E$, let $X_{d} \coloneqq \set{\x \in \monster^{\ell}: d \in
  \acl(E \x)}$. 
$X_{d}$ is \ordefinable (over $E d$) because it is a union of sets defined
by formulae of the form
$\alpha(\cv,d, \x) \et \card{\alpha(\cv, \monster, x)} \leq n $, as $\alpha$ varies among all
$L$-formulae and $n$ varies in $\N$.
Therefore, $X = \bigcup_{d \in \acl(E\av) \setminus E} X_d$ is also \ordefinable (over
$\acl(E \av)$).

\smallskip

By assumption, 
there exists $\bv \in Q$ such that $\rho(\bv/E \av)$ is maximum (inside~$Q$).

Let $\cv \models \tp(\bv/E\av)$ such that $\av \ind_{E \av} \bv$.
As in \cite[(1.10)]{HC:99} (\cf\cite[\Thms 5.12]{MS}), we see that
\[
\cv \ind_{E \bv} \av \et \cv \ind_{E \av} \bv.
\]
Since $\ind$ is \realcanon, and $\acl(E \av) \cap \acl(E \bv) = E$,
we have that
\[
\cv \ind_{E} \av \bv
\]
and therefore $\cv \ind_{E} \av$, proving \ref{prop:HC-claim}.

\medskip

\end{proof}




\section{Simplicity}
\label{sec:simple}

We will now prove that if T is simple, then  $\Tdsg$ is also simple.  Our proof relies on the crucial notion of ``Independence over Algebraically Closed Substructures'' (IACS), which we define and explore in this section.

\begin{thm}[Simplicity of $\Tdsg$]\label{thm:simple} 
If $T$ is simple, then $\Tdsg$ is also simple.
\end{thm}

The above is non-trivial, but thanks to the preliminary work we can give an easy
proof. 

We fix $\tuple{\monster; \deltabar}$ monster model of $\Tdsg$.
Assume that $\monster$ is simple.
Let $\indf$ be Shelah-forking relation on $\monster$,
 and $\indfd$ be the induced independence relation on 
$\tuple{\monster; \deltabar}$ as in \cref{sec:ind}.

Since we want to prove a stronger version of \cref{thm:simple}, we need the following definition.
\begin{definition}\label{def:IACS}
An independence relation $\ind $ satisfies ``Independence over Algebraically
Closed Substructures'' (IACS) if:

\begin{sentence}[(IACS)]
For every $\av, \bv \in \monster^{n}$, $\av', \bv' \in \monster^{m}$,
for every $M \subset \monster$ with $M$ small and algebraically closed, if
\[
\av \ind_{M} \bv, \qquad \av' \ind_{M} \av, \qquad \bv' \ind_{M} \bv, \qquad \av' \elem^{L}_{M}
\bv'
\]
then there exists $\cv \in \monster^{m}$ such that
\[
\cv \ind_{M} \av \bv, \qquad \cv \elem^{L}_{M \av} \av', \qquad \cv \elem^{L}_{M \bv}
\bv'.
\]
\end{sentence}
We say that $T$ satisfies \IACS if $\indf$ satisfies it.
\end{definition}

As we will see later, \IACS is closely related to (weak) elimination of
imaginaries \cf\cite{HC:99}.

\begin{remark}
Obviously, \IACS implies Independence over Models; the converse is not true: for
instance, let $L = \set {E}$, where $E$ is a binary relation, and $T$ be the
theory saying that $E$ has two equivalence classes, both infinite;
then, $T$ is $\omega$-stable and hence simple, therefore $\indf$ satisfies
Independence over Models, but it is easy to see that it does not satisfy \IACS
(take $M = \emptyset$). Indeed it is enough to consider four distinct elements $a, a',
b, b'$ with
$a', b'$ in two different equivalence classes.
\end{remark}

\begin{thm}\label{thm:simple-d}
Assume that $T$ is simple.
Then,
$\indfd$ is a \realstrict independence relation on $\tuple{\monster, \deltab}$ and it
satisfies Independence over Models.

Thus, $\Tdsg$ is simple, and $\indfd$ is Shelah-forking on 
$\tuple{\monster; \deltabar}$.

If moreover $T$ satisfies \IACS, then also $\Tdsg$ satisfies \IACS.
\end{thm}
\begin{proof}
 From \cref{thm:ind} we have that $\indfd$ is a strict independence relation.

We have to show that $\indfd$ satisfies Independence over Models (resp., over
Algebraically Closed Substructures) if $\indf$ satisfies it.
Let $\av$, $\bv$, $\av'$, $\bv'$ be as above and $M$ be a small elementary
$\Ldb$-substructure  (resp., small $\Ldb$-algebraically closed substructure) of $\tuple{\monster, \deltab}$.
The assumptions become:
\begin{align*}
\Jet(\av) &\indf_{M} \Jet(\bv), & \Jet(\av') &\ind_{M} \Jet(\av), \\
\Jet(\bv') &\indf_{M} \Jet(\bv), & \Jet(\av') &\underset{M}{\elem^{\Ldb}} \Jet(\bv')
\end{align*}

By \cref{fact:KP} (applied to~$\monster$), $\indf$ satisfies Independence
over Models (resp., by assumption, $\indf$ satisfies \IACS), and therefore
there exists $\dv \in M^{\Gamma_{n}}$ such
that
\[
\dv \indf_{\monster} \Jet{(\av)} \Jet{(\bv)}, \qquad
\dv \underset{\monster \Jet (\av)}{\elem^{L}}  \Jet{(\av')}, \qquad
\dv \underset{\monster \Jet(\bv)}{\elem^{L}}  \Jet({\bv'})
\]
Let $q \coloneqq \tp^{L}(\dv/ M \Jet(\av) \Jet(\bv))$, and $p$ be restriction of
$q$ to $\monster$.
Notice that $p$ is realized by $\delta \av$, and it is therefore
$\deltab$-compatible.
Moreover, $\dv \indf_{\monster} \Jet({\av}) \Jet({\bv})$ and therefore
$q \indf_{\monster} \Jet(\av) \Jet(\bv)$, and therefore
\[
q \indacl_{\monster} \Jet(\av) \Jet(\bv).
\]
Thus, by the Extension Theorem, $q$ is $\deltab$-compatible, and it is therefore
realized by $\Jet(\cv)$ for some $\cv \in \monster^{m}$.
Thus, we have
\[
\Jet(\cv) \indf_{\monster} \Jet({\av}) \Jet({\bv}), \quad
\Jet(\cv) \underset{\monster \Jet(\av)}{\elem^{L}}  \Jet{(\av')}, \quad
\Jet(\cv) \underset{\monster \Jet(\bv)}{\elem^{L}}  \Jet{(\bv')}
\]
which implies the conclusion.
\end{proof}

\subsection{GEI and simplicity}
\label{sec:simple-GEI}
\begin{thm}\label{thm:simple-GEI}
Assume that $T$ is algebraically bounded, simple, with \GEI and $\indf$ is \realcanon
(\eg $T$ is either stable or supersimple).
Then:
\begin{enumerate}
\item $\indacl = \tind = \indf$;
\item $T$ is super-simple of $SU$-rank 1;
\item $\indacl$ is \eqstrict and \realcanon;
\item Shelah forking on $\Tdsg$ is equal to $\indfd$ and is also \eqstrict and
\realcanon;
\item $\Tdsg$ has \GEI.
\end{enumerate}
\end{thm}
\begin{proof}
$\indacl = \tind$ by \cref{fact:GEI-rosy}
and it is therefore \eqstrict.
$\indacl = \indf$ by \cref{rem:simple-GEI-ind}
$\indacl$ is \realcanon by \cref{cor:indacl-RC}.

$\indfd$ is Shelah's forking by \cref{thm:simple-d}; it is therefore
\eqstrict; it is \realcanon by \cref{thm:ind}.

$\Tdsg$ has \GEI by 3) and \cref{fact:Yoneda}.
\end{proof}
By Theorem~A, we can deduce that the assumptions of \cref{thm:simple-GEI} can be weakened to ``T is algebraically bounded and simple''.



\subsection{Elimination of Imaginaries and simplicity}\label{sec:simple-EI}

In this Section we give some sufficient conditions for $\Tdsg$ to have Elimination of
Imaginaries (EI) when $T$ is simple (see \cref{thm:EI-simple,cor:EI-stable});
for previous results see \cites{MS,Mohamed}.

\medskip

Given a monster model $\monster$, we use the following convention:\\
$\av, \bv, \cv, \dv$ are  finite tuples of real elements,
$A, B,C, D, E$ are small sets of real elements, while $e$ will be an imaginary
element. $\acl$ will always denote the algebraic closure inside $\monster$
(and not $\monster^{eq}$); when we  consider two theories, $T$ and $\Tdsg$,
we will use $\acl$ for the algebraic closure according to $T$, and $\dacl$
for the algebraic closure according to $\Tdsg$.

\begin{definition}\label{good}
We say that $T$ is ``good'' if:
\begin{enumerate}
\item $T$ is algebraically bounded,
\item $T$ has EI,
\item $T$ is supersimple.
\end{enumerate}
\end{definition}

\begin{remark}
If $T$  is good then  $T$ has EHI and $\indf  = \tind$
(\cref{fact:canonical,thm:GEI-rosy-eq}).

Moreover conditions 1) and 2) imply
that $\indacl = \tind$ (\cref{thm:GEI-rosy-eq}).

Conditions 2) and 3) imply that $\indf$ 
satisfies \IACS (\cref{fact:IACS}).

Finally, if $T$ is good, then $T$ is supersimple of SU rank 1.
\end{remark}

\begin{remark}
If we consider the conditions  1) and 2) of Definition \ref{good} together with the assumption that the theory $T$ is stable, then it is easy to see that  $\indf  = \tind = \indacl $ which implies that $T$ is superstable of SU rank 1 (a particular case of supersimplicity), and therefore it is good. 
By Theorem~A, we can further weaken the assumption to conditions 1) and 2) plus $T$ is simple.
\end{remark}


\begin{thm}\label{thm:EI-simple}
If $T$ is good, then $\Tdsg$ is simple,
$\indfd$ is \realcanon and satisfies \IACS,
and
 $\Tdsg$ has EI.
\end{thm}
We will prove the above theorem at the end of this subsection, after some consequences and preliminary results.

\begin{corollary}\label{cor:EI-stable}
If $T$ is algebraically bounded  and stable with EI, then
$\Tdsg$ has EI.
\end{corollary}
\begin{proof}
If $T$ is stable with EI, then it is superstable (of U-rank 1).
(Theorem A gives that $T = \ACFz$).
\end{proof}

\smallskip

We prove a more general result that we will use for the proof of Theorem \ref{thm:EI-simple}.

\begin{proposition}\label{thm:EI-simple-general}
Let $\monster$ be a monster model
satisfying the following conditions:
\begin{enumerate}
\item 
$\monster$ is simple;
\item  $\indf$ satisfies \IACS;
\item $\indf$ is \realcanon;
\item 
$\indf$ is controlled by a rank $\rho$, as defined in 
\ref{en:controlled-rank}
of \cref{prop:HC-claim}.
\end{enumerate}
Then, $\monster$ has Weak Elimination of Imaginaries (WEI).
If moreover $\monster$ expands a field, then $\monster$ has EI.
\end{proposition}

\begin{proof}
The proof keeps following the one in \cite[(1.10)]{HC:99} and \cite[\Thms 5.12]{MS},%
\footnote{Notice that in their proof \cite{MS} do not mention that $\monster$
  should have IACS: but this  holds in their setting.
 \cite[Fact~4.14]{dElbee:23a} gives a similar result that might be easier to
 prove in this context.}.
Here are some more details.

Let $e$ be an imaginary element. Write $e = f(\av)$ for some (real finite tuple)
$\av \in \monster^{\ell}$ and some function $f$ which is definable without parameters.
Let $E := \acleq(e) \cap \monster$.
Our thesis that $U$ has WEI is equivalent to $e \in \dcleq(E)$
(when $\monster$ expands a field, WEI implies EI see Remark \ref{WEIimpliesEI}).

Let $P$ be the set of realization of $\tp(\av/E)$ (inside~$\monster$).
By \cref{prop:HC-claim},
there exists $\cv \in P$ such that $f (\cv) = e$ and $\av$ and $\cv$ are (forking)-independent over~$E$.

\begin{claim}\label{cl:EI-constant-P}
$f$ is constant on $P$.
\end{claim}
Otherwise, there exists $\dv' \in P$ such that $f(\av) \neq f(\dv')$.
Let $\dv'', \av''$ be such that $\dv'' \av'' \indf_{E} \dv' \av$ and
$\dv'' \av'' \elem_{E} \dv' \av$.
Thus, $f(\dv'') \neq f(\av'')$, $\dv'' \elem_{E}  \av$, and
$\av'' \elem_{E} \av$.
Thus, choosing either $\dv \coloneqq \dv''$ or $\dv \coloneqq \av''$ we get
$f(\dv) \neq f(\av) $ and
 and
$\dv \indf_{E} \av$.

By \IACS (since $E$ is algebraically closed), 
(by taking $\av \coloneqq \av$,  $\bv \coloneqq \dv$,
$\av' \coloneqq \cv$, and $\bv' \coloneqq \av$ in \cref{def:IACS})
there exists $\cv' \in \monster^{\ell}$
such that $\cv' \elem_{E \dv} \av$ and $\cv' \elem_{E \av} \cv$.
Since $e = f(\av) \neq f(\dv)$, we get $f(\cv') \neq e$;
but since $e = f(\av) = f(\cv)$, we get $f(\cv') = e$, absurd.

\begin{claim}
$e \in \dcleq(E)$.
\end{claim}
By \cref{cl:EI-constant-P} and compactness, there exists an $L(E)$-definable set
$X$ such that $P \subseteq X$ and $f$ is constant on $X$.
Thus, $e$ is defined by the $L(E)$-formula $\exists \x \in X\, f(\x) = y$.
\end{proof}

\begin{proof}[Proof of \cref{thm:EI-simple}]
We have to show that $T$ and $\Tdsg$ satisfy the assumptions of
\cref{thm:EI-simple-general}.

For $T$, we define  $\rho(\av/B)$ to be the rank given by the algebraic closure
(notice that under our assumptions $T$ is super-simple of rank $1$, and $\rho$ is
equal to the $SU$-rank).
Since $\rho(a_{1}, \dotsc, a_{m}) / B \leq m$, condition  \ref{en:USC}  of Proposition \ref{prop:HC-claim} is automatically
true, and therefore $T$ satisfies the assumptions of \cref{thm:EI-simple-general}.

\smallskip

We consider now $\Tdsg$: we have to show that the various properties of $T$ that
we need are inherited by $\Tdsg$.

We have seen before that $\Tdsg$ is simple with \IACS and \realcanon.

We need to introduce a suitable function $\rho$ and show that it satisfies
the conditions in \cref{thm:EI-simple-general}.
We can use the analogue of the dimension function $\dim_{\mathcal D}$ in
\cite{MS}.
For every $r \in \N$, define $\rho(\av/B)(n)$ as the transcendence degree of
$\Jet_n(\av)$ over $\acl(\Jet(B))$, where 
\[
\Jet_n(\av) = \set{\mu \av: \mu \in \Gamma \mbox{ and } \abs \mu \leq n}
\]
Thus, $\rho(\av/B) \in \omega^{\omega}$.
If we endow $\omega^{\omega}$ with the lexicographic ordering, then $\rho$ is upper
semi-continuous, and satisfies the conditions in
\cref{thm:EI-simple-general}.
\end{proof}


\section{Distality}
This section is based on a suggestion by Elliot Kaplan.\\
Distality is a stronger condition than NIP, it characterizes NIP theories that are as far from stability as possible. For a formal definition and key properties of distal theories, see \cites{Simon,AscheCGZ-22}.
In \cite{FT:24} we showed that if $T$ is NIP, then $\Tdsg$ is NIP.
A similar results hold when $T$ is distal.

\begin{thm}
If $T$ is distal, then $\Tdsg$ is distal.
\end{thm}
\begin{proof}
Immediate from \cref{thm:strong-EQ} and \cite[\Props7.1]{AscheCGZ-22}.
\end{proof}


\section{\texorpdfstring{$\omega$}{\unichar{"1D714}}-stability}\label{sec:w-stable}

In \cite{FT:24}, we showed that if the theory $T$ is stable, then $\Tdsg$ retains
stability. 
The analogous result for $\omega$-stability is false:
here we give necessary and sufficient conditions for $\Tdsg$ to be $\omega$-stable.
Recall that countable theory  is totally transcendental iff the theory is $\omega$-stable (see  \cite[p.125]{Sacks}). 

\begin{thm}\label{thm:w-stable}
$\Tdsg$ is totally transcendental iff $T = \ACFz$ (\ie, $T$ is the theory of pure algebraically closed fields,
possibly with some constants) 
and $\deltab$ commute.
\end{thm}

We first need the following lemma (\cf Theorem~A):
\begin{lemma}
$T$ is algebraically bounded and strongly minimal iff $T = \ACFz$.
\end{lemma}
\begin{proof}
It is clear that, if $T = \ACFz$, then $T$ is algebraically bounded and strongly minimal.

Conversely, assume that  $T$ is algebraically bounded and strongly minimal.
Let $\K \models T$.
Since $\K$ is strongly minimal, as a field $K$ is algebraically closed (see
Macintyre's Theorem \cite[Thm.3.1]{Poizat:groups}).
By \cite[Thm.1]{Hrushovski:92}, every $\K$-definable set is definable in the field language.
\end{proof}





\begin{proof}[Proof of \cref{thm:w-stable}]
If $T = \ACFz$, then $\Tdbg$ is $\omega$-stable  \cite{mcgrail}.

Let us now consider the converse.
Let $\Kdb \models \Tdsg$.

\begin{enumerate}
\item Assume that $\deltab$ do not commute.
Let $C_{1}$ be the fixed field of $\delta_{1}$.
Notice that $\delta_{2}(C_{1}) = K$, since $\K$ is existentially closed. 
Thus, $\MoR(C_{1}) = \MoR(K)$.
However, $C_{1}$ is a subgroup of $K$ of infinite index, and therefore 
 $\MoR(K) = \infty$, proving that $\Kdb$ is not $\omega$-stable.

\item We have to prove that $T = \ACFz,$ so by Lemma  \ref{thm:w-stable} we assume that $K$ is not strongly minimal.
We want to show that $\Kdb$ is not $\omega$-stable.
It suffices to show that, assuming that $L(K)$  is countable,
there exist $2^{\omega}$ $\Ldb$-1-types over~$K$.
Then, let $X \subset K$ which is $\K$-definable and such that both
$X$ and $\K \setminus X$ are infinite.
Fix $J \subseteq \N$ and
consider the following partial $\Ldb(K)$-type.
\[
p_{J}(x) = \Pa{\delta_{1}^{i} x \in X: i \in J} \cup \Pa{\delta_{1}^{i} x \notin X: i \notin J}.
\] 
Notice  $p_{J}$ is indeed a partial type and, for every
$J \neq J'$, $p_{J}$ and $p_{J'}$ are incompatible.
Thus, $\card{S^{1}_{\Ldb}(K)} \geq 2^{\omega}$.
\qedhere
\end{enumerate}
\end{proof}


\section{The field of constants}\label{sec:constant}

We give some interesting results on the field of constants. Before we introduce some notations and definitions.


\begin{definition}
The \intro{field of constants} is the set
\[
\cfieldK \coloneqq \set{a \in K: \delta_{1}(a) =  \dots = \delta_{k}(a) = 0}
\]
\end{definition}

We consider the reduct $\tuple{\K,\cfieldK}$ (that is, the
expansion of $\K$ with a unary predicate for $\cfieldK$).

Observe that $\cfieldK$ is \intro{dense} in $\K$ \wrt the matroid $\acl$: that
is, for every $Z \subseteq K$ which is $L$-definable (remember that $L$ is the signature of~$\K$) with parameters in $\K$ and large,
$Z$ intersects $\cfieldK$; moreover,
$\cfieldK$ is also L-algebraically closed in~$\K$.
Thus, $\tuple{\K,\cfieldK}$ is a \textbf{lovely pair of
geometric structures} (in the sense of \cites{BerensteinV:10,Boxall}: see also
\cite{Fornasiero:matroid}).
Thus, we can apply the known results (see \cites{BerensteinV:10,Boxall, Fornasiero:matroid}).

\begin{lemma}\label{lem:P-indep}
    $\bv \in K^{\ell}$ and $B := \dcl_L(\Jet(\bv)= \dcl_{\Ldb}(\Jet(\bv))$. 
    Then, $B$ is P-independent: that is, $B \indacl_A \cfieldK$, where
    $A := B \cap P$.
\end{lemma}
\begin{proof}
Assume not: then, there exists $\dv 
\in \Jet(\bv)^n$ 
s.t.\ $\dv$ is L-algebraically independent over $A$ but
the following system has a solution $\cv \in (K^n)^m$:
\begin{equation}\label{eq:d-system}
\left\{\begin{aligned}
    \delta \x &= 0\\
    \x &\neq 0\\
    \sum_{i,j} x_{ij}d_j^i &= 0.
\end{aligned}\right.
\end{equation}
\Wlog, we may assume that $\tuple{\K, \deltab}$ is a monster model.
Let $\bv' \elem_A^{\Ldeltabar} \bv$ s.t. $\bv' \indacld_A \cv$.
Let $\sigma$ be an $\Ldb$-automorphism of $\K$ over $A$ mapping $\bv$ to $\bv'$ and $\dv' := \sigma(\dv)$.
Then, $\dv'$ is $L$-algebraically independent over $A$ and therefore
\eqref{eq:d-system} (with $\dv'$ in place of $\dv$) 
cannot have a solution in~$K$.
But $\bv$ and $\bv'$ have the same $\Ldb$-type, absurd.
 \end{proof}

\begin{thm}
$\cfieldK$ is an elementary $L$-substructure of $\K$.

Let $\bv \in K^{\ell}$  and $X \subseteq \cfieldK^{n}$ be $\Ldb$-definable with parameters~$\bv$.
Then, there exists $Y \subseteq K^{n}$ which is $L$-definable (in $\K$) 
with parameters $\Jet(\bv)$
such that $X = Y \cap \cfieldK^{n}$.
If moreover $\bv \in \cfieldK$, then $X$ is $L$-definable in $\cfieldK$ with
parameters $\bv$: equivalently, there exists $Y \subseteq K^{n}$ which is $L$-definable
in $\K$ with parameters $\bv$ such that $X = Y \cap K^{n}$.
\end{thm}
\begin{proof}
The Theorem is a reformulation of \cite[Lemma 8.22]{Fornasiero:matroid} (see also \cite[\Props 3.4]{BerensteinV:10}), using \cref{lem:P-indep}.
\end{proof}
Notice that, by \cref{cor:acld}, $\cfieldK$ is algebraically closed in $\K$ \wrt the $\Ldb$-structure.


\begin{definition}
A basic formula is a formula of the form
\[
\exists \y\, \Pa{\y \in K^{\ell} \wedge \psi(\x, \y)}
\]
where $\psi$ is an $L$-formula.
A basic set is a set definable by a basic formula (with parameters from $\K$).
\end{definition}

\begin{thm}
Let $Z \subseteq K^{n}$ be definable in $\tuple{\K,\cfieldK}$ with parameters
from~$\K$.
Then, $Z$ is a finite Boolean combination of basic sets, with the same
parameters as $Z$.
\end{thm}
The Theorem is a reformulation of \cite[\Cors 3.2]{BerensteinV:10} and \cite[\Thms 8.5]{Fornasiero:matroid}.

\begin{remark}
Let $\tuple{A,B}$ be a lovely pair of geometric structures, with $A \models T$.
Then, there exists $B^{*} \succeq B$ and a derivation $\delta^*$ on $B^{*}$ such that
$\tuple{B^{*}, \delta^{*}} \models \Tdg$ and $\tuple{B^{*}, A^{*}} \succeq \tuple{B,A}$, 
where $A^{*} \coloneqq \mathfrak C_{\delta^{*}} $ 
\end{remark}
\begin{proof}
The theory $T^{lovely}$ of lovely pairs of models of $T$ is complete (see \cite{BerensteinV:10}).
Let $\tuple{B^{*}, A^{*}} \succeq \tuple{B,A}$ be a 0-big 
(a.k.a. ``splendid'': see \cite{Hodges}).
By bigness, there exists a derivation $\delta^{*}$ on $B^{*}$ satisfying the conclusion.
\end{proof}
\cite[\S5]{KP} use a particular case of the above remark to (re-)prove some
results about lovely pairs.

For more results (in particular on imaginaries in $\tuple{\K,\cfieldK}$) see
\cites{BerensteinV:10, Boxall, Fornasiero:matroid}.

\section{Open core}\label{sec:open}

The open core of a topological structure  is the structure generated by the sets which are
both open and definable.
We investigate expansions of algebraically bounded topological structures with generic derivations:
more precisely we prove that, under very weak assumptions, $T$~is the open core
of $\Tdsg$ (\cref{cor:OC}).

By a topology $\tau$ on  $T$  we mean an assignment, to every
$\K \models T$ and every $n \in \N$, of a topology $\tau_{n, \K}$ on $K^{n}$.



\begin{definition}[{\cite{Pillay:87}}]
$\tau$ is definable if, for every $n\in \N$, there exists an $L$-formula
(without parameters) $\alpha_{n}(\x,\y)$, where $\x$ is an $n$-tuple and
$\y$ is a finite tuple,
such that, for every $\K \models T$, we have that
\[
\set {\alpha(\K,\bv): \bv \in K^{\abs{\y}}} 
\]
is a basis of open sets for $\tau_{n, \K}$.
The family $\set {\alpha(\K,\bv): \bv \in K^{\abs{\y}}}$ is a definable basis for the topology.
\end{definition}

If we give only a topology on $K$, we assume that the topology on $K^{n}$ is
the product topology.

\begin{definition}
$T$ is an open core of $\Tdsg$ if, for every $\Kdb \models \Tdsg$ and every $n \in \N$,
every $\tau$-open $\Ldb(K)$-definable
subset of $K^{n}$ is already $L(K)$-definable.  
\end{definition}
We give: sufficient conditions for $T$ being an open core of $\Tdsg$
(see \cref{ass:open-core,thm:open-core}) 
and for
$\Tdsg$ to have Elimination of Imaginaries (see \cref{thm:EI}).

In an article in preparation we will present another approach,
with slightly different assumptions and a very
different proof.


\bigskip

Let $\K \models T$ and $n \in \N$.
Let $X \subseteq K^{n}$ be an $L(K)$-definable set.

\begin{definition}
Given $\av \in K^{n}$, 
the \intro{Local Dimension} of $X$ at $\av$ is 
\[\ldim_{\av}(X) \coloneqq \min\set{ \dim(X \cap U): U \text{ open definable set containing $\av$}}.
\]
The Local Dimension of $X$ is
\[
\ldim(X) \coloneqq \max\set{ \ldim_{\av}(X): \av \in X }.
\]
We say that $X$ has \intro{constant local dimension}~$d$ if, for every
$\av \in X$, $\ldim_{\av}(X) = d$.
\end{definition}
\begin{remark}
$\ldim(X) \leq \dim(X)$.
\end{remark}

\begin{definition}
$\dim$ is local (\wrt to $\tau$) if, for every 
$\K \models T$, for every $n \in \N$, and every $X \subseteq K^{n}$ $L(K)$-definable,
 $\dim(X) = \ldim(X)$.
\end{definition}

\cite{FH:12} gives sufficient conditions for dim to be local. 
It gives also an example where it is not local (the Sorgenfrey plane).
\begin{fact}[\cite{FH:12}]\label{fact:dim-local}
Assume that $\tau$ is definable, the topology on $K^{n}$ is the product topology
(for every $n \in \N$),
$\tuple{K,+}$ is a Hausdorff topological group (for every $\K \models T$),
and every open definable set is large, i.e. it has the same dimension of ambient space.

Then, $\dim$ is local.
\end{fact}

\begin{assumptions}\label{ass:open-core}
We assume that, for every $n,m \in \N$ and $\K \models T$,
\begin{enumerate}
\item $\tau$ is definable;
\item the projection map $K^{n} \times K^{m} \to K^{n}$ is continuous;
\item $\tau_{n, \K}$ is invariant under permutation of coordinates;
\item $\dim$ is local.
\end{enumerate}
\end{assumptions}

\begin{examples}
The following topologies satisfy the assumptions.
\begin{enumerate}
\item $T =\RCF$ with the usual Euclidean topology.
\item $T$ is a theory of Henselian valued fields with the valuation topology.
\end{enumerate}
\end{examples}

\begin{thm}[Open Core]\label{thm:open-core}
$T$ is the open core of $\Tdsg$.

More precisely: let $\Kdb \models \Tdsg$ and $A \subseteq K$ s.t. $\deltab(A) \subseteq A$.
Let $X \subseteq K^{n}$ be $\Ldb(A)$-definable and closed.
Then, $X$ is $L(A)$-definable.
\end{thm}

We give a proof later in Section \ref{proofEI}. Particular cases of \cref{thm:open-core} were already known: 
see \cite[\S6]{KP}.

\begin{thm}\label{thm:EI}
Assume  moreover that $\tau$ 
satisfies the following condition:
\begin{enumerate}[start=5]
\item\label{top:boundary} If $X$ is $L(\K)$-definable and nonempty, then $\dim(\overline X \setminus X) <
\dim (X)$, where $\overline X$ is the topological closure of $X$.
\end{enumerate}
Then, $\Tdsg$ has Elimination of Imaginaries modulo $T^{eq}$.
\end{thm}
\begin{proof}
With trivial modifications, the proof suggested by M. Tressl in \cite{FK} works.
\end{proof}

For the following corollary, we spell all conditions explicitly.
\begin{corollary}\label{cor:OC}
Let $T$ be algebraically bounded, with a definable topology~$\tau$.
Assume that
 $\tau$ is a non-trivial and non-discrete ring topology.
Then, $T$ is the open core of $\Tdsg$.

If moreover $\tau$  and satisfies  Assumption~\ref{top:boundary} in \cref{thm:EI}, then
$\Tdsg$ has Elimination of Imaginaries modulo $T^{eq}$.
\end{corollary}
\begin{proof}
Since $\tau$ is a non-trivial ring topology, it is Hausdorff (see \eg \cite[\Thms 11.10]{Warner:89}).

We only need to show that $\dim$ is $\tau$-local.

By \cref{fact:dim-local}, it suffices to show that every nonempty open  subset of $K^{n}$ is large.

For $n=1$, since $\tau$ is Hausdorff and non-discrete, every nonempty open subset
of $K$ is infinite and hence large.

For $n > 1$, let $U \subseteq K^{n}$ be open and nonempty.
Let $B_{1}, \dotsc, B_{n} \subseteq K$ be open and nonempty s.t. $B \coloneqq B_{1} \times \dots \times B_{n}
\subseteq K^{n}$.
Since, by the case $n=1$, each $B_{i}$ is large, then $B$ is large.
\end{proof}

\subsection{Proof of Theorem \ref{thm:open-core}}\label{proofEI} We present some preliminary results before giving the proof.

Fix $n \in \N$ and $\x = \tuple{x_{1}, \dotsc, x_{n}}$.
We use the same notation as in \cref{notazioni}.

We need to endow $K^{\Gamma_{n}}$ with a topology.
For every $J$ finite subset of $\Gamma_{n}$ of cardinality $m \in \N$, $\tau_{m, \K}$ is a
topology on $K^{m}$ and hence on $K^{J}$ (it does not depend on how we enumerate the elements of~$J$,
since we assumed that $\tau_{m, \K}$ is invariant under permutation of coordinates).
We have a natural map $\Pi_{J}: K^{\Gamma_{n}} \to K^{J}$.
We define the topology $\tau_{\K,\Gamma_{n}}$ as the coarsest topology making all the maps
$\Pi_{J}$ continuous: equivalently, a basis of $\tau_{\K, \Gamma_{n}}$ is given by
\[
\set{\Pi_{J}^{-1}(U): U \subseteq K^{J} \text{ open}, J \subset_{fin} \Gamma_{n}}.
\]

Let $\K \models T$ be a monster model (\ie, $\lambda$-saturated and $\lambda$-homogeneous for some
sufficiently large cardinal $\lambda > \aleph_{0} + \card L$) and $A \subset K$ be a ``small''
subset (\ie, $\card A < \lambda$).

Let $J \subseteq \Gamma_{n}$ ($J$ could be finite or infinite);  again, we consider $J$ as a
set of indices and $x_{J}$  as a 
tuple of variables (see Section \ref{notation} )

A partial type $p(x_{J})$ over $A$ determines a (type-definable) subset 
$Z \subseteq K^{J}$.
If $J$ is finite,
we can define the dimension of $Z$ as
\[
\dim(Z) \coloneqq \min \set{\dim(\alpha(\K): \alpha \in p} \in \N.
\]
If however $J$ is infinite, the dimension of $Z$ (defined in the ``obvious'' way)
may be infinite (e.g. $\dim(K^{J})$ is infinite).
We need to define when $Z'$ is a large subset of $Z$ even when $\dim(Z)$ is
infinite.
\begin{definition}
Let $Y \subseteq Z$ be subsets of $K^{J}$ which are type-definable (over some small 
$A \subset K$).
We say that $Y$ is \intro{large in} $Z$ if, for every $J' \subseteq J$ with $J'$ finite,
$\dim(\Pi_{J'}(Y)) = \dim(\Pi_{J'}(Z))$.
\end{definition}

Given a partial type $p(x_J)$, we denote by $p(\K)$ 
the set of realizations of $p$ in $\K^{J}$, and
$\dim(p(x_J)) \coloneqq \dim(p(\K)$).
Given $J' \subseteq J$, we denote by $p(J')$ the partial type corresponding to the
set $\Pi_{J'}(p(\K))$.

Let $q \in S^{J}(A)$ be a complete type, of dimension~$d$.

\begin{lemma}\label{lem:const-dim}
Assume that 
$J$ is finite.

Then, there is a family $\set{X_{i}: i \in I}$ of $L(A)$-definable sets such that:
\begin{enumerate}
\item $q(\K) = \bigcap_{i \in I} X_{i}$;
\item each $X_{i}$ has constant local dimension~$d$.
\end{enumerate}
\end{lemma}
\begin{proof}
Let $\set{Y_{i}:  i \in I}$ be a family of $L(A)$-definable sets such that:
$q(\K) =  \cap \set{Y_{i}: i \in I}$.
Let $Z$ be some $L(A)$-definable set of dimension $d$ containing $q(\K)$: by
replacing $Y_{i}$ with $Y_{i} \cap Z$, \wloG we may assume that all $Y_{i}$ have
dimension~$d$.
Given an $L(A)$-definable set $Y \subseteq K^{J}$ of dimension~$d$, let
\[
\sing(Y) \coloneqq \set{\av \in Y: \ldim_{\av}(Y) < d}\\ \mbox{ and } 
\reg(Y) \coloneqq Y \setminus \sing(Y).
\]
Since $\dim$ is definable, $\sing(Y)$ and $\reg(Y)$ are $L(A)$-definable.
Since $\dim$ is local, $\dim(\sing(Y)) < d$,
and therefore $\dim(\reg(Y)) = d$ and $\reg(Y)$ has constant local
dimension~$d$.
\begin{claim}
For every $i \in I$, $q(\K) \subseteq \reg(Y_{i})$.
\end{claim}
Assume not: then, since $q$ is a complete type over $A$,
$q(\K) \subseteq \sing(Y_{i})$, but then $\dim(q) \leq \dim(\sing(Y_{i})) < d$, absurd.

Thus, we can define $X_{i} \coloneqq \reg(Y_{i})$.
\end{proof}
Notice that in the above lemma we used that $\dim$ is definable.

\begin{lemma}\label{lem:local-dim-type}
Let $U \subseteq K^{J}$ be open and definable (with parameters in some small set $B$ with
$A \subseteq B \subset K$).
If $U \cap q(\K^{J})$ is nonempty, then it is large inside $q(\K^{J})$.
\end{lemma}
\begin{proof}
Since $U$ is definable, there is some $J_{0} \subseteq J$ finite and 
$V \subseteq K^{J_{0}}$ $L(A)$-definable and open, such that
$U = \pi_{J_{0}}^{-1}(V)$.
Let $J_{0} \subseteq J' \subseteq J$, with $J'$ finite.
We need to prove that $\Pi_{J'}(U \cap q(K^{J'}))$ is large inside 
$\Pi_{J'}(q(\K))$.

The former is equal to $U' \cap q(\K(J'))$, where
$U' \coloneqq \Pi_{J'}(U')$: thus, by replacing $U$ with $U'$,  $J$ with $J'$, and
$q(x_{J})$ with $q(x_{J'})$, we may assume that $J'$ is finite.
Thus, by \cref{lem:const-dim}, we may assume that $q(\K) = \bigcap_{i}X_{i}$,
where each $X_{i}$ is an $L(A)$-definable set of constant local dimension~$d$.
Thus, since $X_{i} \cap U$ is nonempty, $X_{i} \cap U$ has (local) dimension $d$, and
therefore
$\dim(\bigcap _{i} X_{i} \cap U) = d$.
\end{proof}


\begin{proof}[Proof of Theorem \ref{thm:open-core}]
\Wlog, we may assume that $\Kdb$ is a monster model and $A$ is ``small''.
By Beth definability, it suffices to show that $X$ is $L(A)$-invariant
(\ie, that $X$ is set-wise invariant under automorphisms of $\K$ as $L$-structure fixing $A$
point-wise).

Let $Z \coloneqq \Jet(X) \subseteq K^{\Gamma_{n}}$ and $Y$ be the closure of $Z$ (according to
the topology $\Pi_{\K, \Gamma_{n}}$).
Since $X$ is closed, we have that $X = \Pi_{n}(Y)$ (where $\Pi_{n}$ is the
projection onto the first $n$ coordinates).
Thus, it suffices to show that $Y$ is $L(A)$-invariant.

Let $\av, \av' \in K^{\Gamma_{n}}$, with  the same $L(A)$-type, and such that
$\av \in Y$.
We need to show that $\av' \in Y$.
Since $\tau_{\K, n}$ is definable and $Y$ is closed, it suffices to prove the
following:
\begin{claim}
Let $U' \subseteq K^{\Gamma_{n}}$ be an $L(\K)$-definable open set containing~$\av'$.
Then, $U'$ intersects~$Y$.
\end{claim}
Let $\phi$ be an $L(A)$-automorphism of $\K$ such that $\phi(\av') = \av$.
Let $U \coloneqq \phi(U')$: notice that $U$ is an open $L(A)$-definable set
containing $\av$; therefore, there exists $\bv \in X$ such that
$\Jet(\bv) \in U$.
Let $\cv'\coloneqq \phi^{-1}(\Jet(\bv))$.

Let $q(\x_{\Gamma})$ be the $L(A)$-type of $\Jet(\bv)$: by definition, it is also the
$L(A)$-type of $\cv'$.
By \cref{lem:local-dim-type}, $q(\K) \cap U'$ is large inside $q(\K)$.
Thus, by the Extension Theorem, there exists $\bv' \in K^{n}$ such that
\[
\Jet(\bv') \in q(\K) \cap U'.
\]
Therefore, $\bv' \in X$ and therefore
$\Jet(\bv') \in Z \cap U' \subseteq Y \cap U'$.
\end{proof}

\section{Differential dimension}\label{sec:ddim}

\subsection{The commutative case}\label{sec:ddim-comm}

Let $\tuple{A, \deltabar}$ be a field (of characteristic~$0$) 
with $k$ commuting derivations.
The derivations induce a matroid on $A$.
Given $a \in A$, $Y \subseteq A$ and  $X\subseteq A$, we define $a \in \deltacl_{Y}(X)$ if $\Jet(a)$ is not
algebraically independent over $\Jet(X \cup Y)$.

As shown in \cite{FK}, $\deltacl_{Y}$ is a matroid on $A$.%
\footnote{A more general result is true. Let $\mathbb A := \tuple{A, \cl}$ be a finitary
  matroid. Let $\deltabar$ be a tuple of commuting quasi-endomorphisms of
  $\mathbb A$, in the sense of \cite{FK}. Given $X, Y \subseteq A$, define
$\deltabar$-$cl_{Y}(X)$ as the set of $a \in A$ such that $a^{\Jet}$ is not $\cl$-independent
over $\Jet (X)\Jet(Y)$. Then, $\deltabar$-$\cl$ is a finitary matroid  on~$A$.
}

Fix  $\tuple{\monster, \deltabar}$ monster model of $\Tdbg$.
We have the corresponding matroid $\deltacl(X) \coloneqq \deltacl_{F}(X)$.
\begin{thm}
$\deltacl$ is an existential matroid (in the sense of \cite{Fornasiero:matroid}).
\end{thm}
\begin{proof}
We have to prove that $\deltacl$ is definable and it satisfies existence
(see \cite[\S3]{Fornasiero:matroid}).

The fact that $\deltacl$ is definable  means that, for every $A \subseteq \monster$ and
$b \in \deltacl(A)$
there exists an $\Ldeltabar$-formula $\phi(\x, y)$ and $\av \in A^{n}$
such that $\tuple{\monster, \deltabar} \models \phi(\av,b)$ and, for every 
$\av',b'$ in $\monster$, if $\tuple{\monster, \deltabar} \models \phi(\av', b')$, then
$b' \in \deltacl(\av')$.
We can take as $\phi$ any formula witnessing that $b^{\Jet}$ is not algebraically
independent over~$A$.

\smallskip

For existence, let $A \subseteq B \subset \monster$ be subsets of small cardinality.
Let $c \in \monster$ such that $c \notin \deltacl(A)$.
We have to show that there exists $d \in \monster$ such that $c$ and $d$ have the
same $\Ldelta$-type over $A$ and $d \notin \deltacl(B)$.

Since $\indacld$ satisfies existence, 
there exists $d \in \monster$ such that $c$ and $d$ have the same 
$\Ldelta$-type over $A$ and $d \indacld_{A} B$.
Then, $\Jet(d)$ is algebraically independent over $\Jet(B)$: therefore, 
$d \notin \deltacl(B)$, proving that $\deltacl$ is an existential matroid.
\end{proof}
Thus, $\deltacl$ induces a dimension function $\deltadim$ on models of $\Tdbg$
(see \cite{Fornasiero:matroid}; see also \cite{GP:12}).

\begin{remark}
\begin{enumerate}
\item  $\deltacl$ is not the $\Tdbg$-algebraic closure: the former only
contains the latter. 
For instance, the whole field of constants $\cfieldK$ is in
$\deltacl(\emptyset)$.
\item   $\indacld$ is not the independence relation induced
by~$\deltacl$, because $\indacld$ is strict. 
For instance, if
$a \in \cfieldK \setminus \acl(F)$, then $a \notind[$M,\deltab$]_{\emptyset} a$.
\end{enumerate}
\end{remark}

\begin{lemma}[See \cite{ELR}]
Let $\tuple{\K, \deltabar} \models \Tdbg$.
Let $Y \subseteq K^{n}$ be $L$-definable (with parameters).
Then, $\dim(Y) = \deltadim(Y)$.
\end{lemma}
\begin{proof}
By the properties of dimension functions (see \cite{Dries:89}) it suffices to
treat the case when $n = 1$ (the general case follows by induction on $n$).
If $\dim(X) = 0$, then $X$ is finite, and therefore $\deltadim(X) = 0$.
If $\dim(X) = 1$, then $(X - X)/(X - X) = \K$, and therefore $\deltadim(X) = 1$.
\end{proof}

The same proof gives a more general result.
\begin{proposition}[Invariance of dimension for fields]
Let $L$ be a language expanding the language of rings, and
$L^{*}$ be an expansion of~$L$.
Let $A^{*}$ be an $L^{*}$-structure expanding a field, and $A$ be its
restriction to the language~$L$.
Assume that $\dim^{*}$ and $\dim$ be dimension functions on~$A^{*}$ and~$A$,
respectively.
Then, for every $X \subseteq A^{n}$ which is $L$-definable (with parameters),
$\dim^{*}(X) = \dim(X)$.
\end{proposition}

\smallskip

Unlike in the case of lovely pairs, we cannot approximate $\Ldelta$-definable
sets with $L$-definable sets.
\begin{remark}
Let $\tuple{\K, \delta} \models \Tdg$.
Let $X \subseteq K$ be $\Ldelta$-definable (with parameters).
If $X$ is definable in the lovely pair $\tuple{\K, \cfieldK}$ (see \cref{sec:constant}),
then there exists $Y \subseteq K$ which is $L$-definable and such that $\deltadim(X \Delta
Y) < 1$ (\cite[Proposition 8.36]{Fornasiero:matroid}).
If not, such $Y$ might not exist: for instance, let $\K$ be a real closed field,
and $X \coloneqq \set{x \in K: \delta x > 0}$.
\end{remark}

\subsection{The non-commutative case}

The assumption that the derivations commute cannot be dropped.
\begin{lemma}
If the derivations do not commute, then $\deltacl_{Y}$ is not a matroid, because it
is not transitive.%
\footnote{Naturally, we use the free monoid $\Gamma$ instead of the free commutative
  monoid $\Jet$ to define $\deltacl$ in this situation.}
In fact, let $k = 2$ and $\tuple{\K, \deltabar} \models \Tdng$.
Then, there exist $a, b, c \in K$ such that:
\begin{enumerate}
\item $a^{\Gamma}$ is algebraically independent over $F$;
\item $\delta_{2} b = 0$ and $\delta_{1} b = \delta_{1} a$; 
\item $c = a - b$.
\end{enumerate}
Notice that $\delta_{1}c = 0$.
Then, $a \notin \deltacl(F)$, $b,c \in \deltacl(F)$, but $a \in \deltacl_{F}(b,c)$: thus,
transitivity fails.
\end{lemma}

\begin{lemma}
For $ k\geq 2$, models of $\Tdng$ do not have a dimension function.
\end{lemma}
\begin{proof}
For simplicity, we do the case when $k = 2$.
Let $\tuple{\K, \deltabar} \models \Tdng$.
Assume, by contradiction, that $\dim'$ is a dimension function on
$\tuple{\K, \deltabar}$. 
Let $X \coloneqq \set{b \in K: \delta_{1} b = 0}$.
Let $Y \coloneqq \set{c \in K: \delta_{2} c = 0}$.
Notice that $X$ and $Y$ are $\Ldelta$-definable subfields of $K$ of infinite
index inside $\K$: thus, $\dim'(X) = \dim'(Y) = 0$.
\begin{claim}
$X + Y = K$.
\end{claim}
Let $a \in K$.
Let $b \in K$ such that $\delta_{1} b = 0$ and $\delta_{2} b = \delta_{2} a$, and let $c = a - b$.
Notice that $b \in X$ and $c \in Y$.
Thus $X + Y = K$.

But $\dim'(X) = \dim'(Y) = 0$, and therefore $\dim'(K) = 0$, while the axioms
of dimension require that $\dim'(K) = 1$.
\end{proof}

\section{Genericity}\label{sec:Polish}

We denote by $K^{K}$  the set of all functions from $K$ to $K$, and by
$\DerK \subset K^{K}$ the set of derivations on $K$ extending~$\eta$
(remember that $\eta$  is a fixed derivation
on the field $F = dcl(\emptyset)$: see \S\ref{subsec:prelim-derivation}).
The main references for this section is  \cite{Hjorth}, from which our presentation is heavily inspired; for the background notions of descriptive set theory see \cite{Kechris}.

For every $\av, \bv \in K^{n}$, we define
\[
B_{\av, \bv} \coloneqq \set{\delta \in K^{K}: \delta(\av) = \bv}.
\]
For every $\Ldelta$-sentence $\phi$ with parameters in $\K$, we define
\[
U_{\phi} \coloneqq \set{\delta \in K^{K}: \tuple{\K, \delta} \models \phi}.
\]
The set $K^{K}$
has two ``canonical'' topologies:
\begin{itemize}
\item  The pro-discrete topology, whose basis of open sets is given by
\[
\set{B_{\av,\bv}: \av, \bv \in K^{n}, n \in \N},
\]
and which we denote by $\taud$, which is the topology induced by the product topology on $K^{K}.$
\item
The  ``first-order'' topology, whose basis of open sets is given by
\[
\set{U_{\phi}: \phi \text { $\Ldelta$-sentence with parameters in $\K$}},
\]
and which we denote by $\tauFO$.
\end{itemize}
\begin{remark}
Another basis for $\taud$ is
\[
\set{U_{\phi}: \phi \text { quantifier-free
    $\Ldelta$-sentence with parameters in $\K$}}
\]
\end{remark}
In fact,
\[
  B_{\av,\bv} = U_{(\delta a_{1} = b_{1} \wedge \dotsb \wedge \delta a_{n} = b_{n})}.
\]

For the remainder of this section, when we don't specify the topology, we mean~$\taud$.
The following notions concerning specific formulas are useful for obtaining information about sets with respect to the topology.   

\smallskip

\begin{definition}
We say that an $\Ldelta$-sentence $\phi$ with parmeters $\av$
is ``relatively quantifier free'' if
$\phi = \alpha(\Jetd{}(\av))$ for some $L$-formula without parameters~$\alpha$.\\
We say $\phi$ is ``relatively existential''   if 
\[
\phi = \exists \x\ \alpha(\Jet(\av), \Jet(\x)),
\]
for some $L$-formula without parameters $\alpha$.\\
We can define
``relatively universal'' and ``relatively $\forall\exists$'' $\Ldelta$-sentences with
parametrs.

\end{definition}

\begin{lemma}\label{lem:topology-formulae}
Let $\phi$ be an $\Ldelta$ sentence with parameters.
\begin{itemize}
\item If $\phi$ is relatively quantifier free, then $U_{\phi}$ is clopen.
\item If $\phi$ is relatively  existential, then $U_{\phi}$ is open.
\item If $\phi$ is relatively  universal, then $U_{\phi}$ is closed.
\item If $\phi$ is relatively  $\forall\exists$, then $U_{\phi}$ is $\Gdelta$.
\end{itemize}
\end{lemma}
\begin{proof}
We do only the case when $\phi$ is relatively existential: the others are similar.
Write $\phi = \exists \y\ \alpha(\av, \delta \av, \dotsc, \delta^{n}\av, \bv, \delta \bv, \dotsc, \delta^{m} \bv)$,
for some $L$-formula $\phi$.
Then,
\begin{multline*}
U_{\phi} = \bigcup \Pa{B_{\av,\av_{1}} \cap  B_{\av_{1},\av_{2}} \cap \dotsb \cap
B_{\av_{n-1},\av_{n}} \cap B_{\bv,\bv_{1}} \cap B_{\bv_{1},\bv_{2}} \cap \dotsb \cap
B_{\bv_{m-1},\bv_{m}}:\\ 
\av_{1}, \dotsc, \av_{n}, \bv_{1}, \dotsc, \bv_{m} \in K^{< \omega} \wedge 
\tuple{\K, \delta} \models \alpha(\av, \av_{1}, \dotsc, \av_{n}, \bv, \bv_{1},
\dotsc, \bv_{m})}. 
\end{multline*}
\end{proof}

For the remainder of this section, we assume that $\K$ and $L$ are \textbf{countable}.

\smallskip

Thus, $\DerK$ is $\tau_{FO}$-closed and it is a $\taud$-$\Gdelta$ inside $K^{K}$;
we use the same names for the induced topologies on~$\DerK$.
Notice that $K^{K}$  is a Polish space: therefore, $\DerK$
is also a Polish space (see \cite{Kechris}).
Thus, any two dense $\Gdelta$ subsets of $\DerK$ always intersect.

Given $Z \subseteq K^{n} \times K^{n}$, we define
\[
I_{Z} \coloneqq \set{ \delta \in \DerK: \exists \bv \in K^{n}: \tuple{\bv,\delta \bv} \in Z }
\]

\begin{lemma}\label{open}
For every $Z \subseteq K^{n} \times K^{n}$, $I_{Z}$ is an open subset of $\DerK$.
\end{lemma}
\begin{proof}
\[
I_{Z} = \bigcup\Pa{B_{\av,\bv} : \av \in K^{n}, \bv \in K^{n}, \tuple{\av, \bv} \in Z}.
\]
\end{proof}

Let $\mathbb G$ be the family of derivations $\delta \in \DerK$ such that 
$\tuple{\K, \delta} \models \Tdg$.
Let $\mathcal L$ be the family of the sets $Z \subseteq K^{n + n}$ definable with
parameters, such that $\Pi_{n}(Z)$ is large (for some $n \in \N$).  
\begin{lemma}
$\mathbb G(M) = \bigcap_{Z \in \mathcal L} I_{Z}$.
Moreover, $\mathbb G$ is a $\Gdelta$-subset of $\DerK$.
\end{lemma}
\begin{proof}
By the axiomatization $\Tdgpp$ (where we introduced in \cite[Def.~3.3]{FT:24} ), $\mathbb G(M) = \bigcap_{Z \in \mathcal L} I_{Z}$.
Each $I_{Z}$ is open.
By our assumptions, $\mathcal L$ is countable.
\end{proof}

\begin{lemma}
On $\Generic$, $\tauFO$ and $\taud$ coincide.
\end{lemma}
\begin{proof}
By elimination of quantifiers, every $\Ldelta$-sentence is equivalent, modulo
$\Tdg$, to a relatively quantifier-free sentence.
The conclusion follows from \cref{lem:topology-formulae}. 
\end{proof}

\begin{lemma}
Assume that $\rk(K / F)$ is infinite.
Then, for every $\av$ finite tuple in $K$ and every $W$  large subset of
$K^{n}$ which is $L$-definable with parameters,
there exists $\bv \in W$ which is algebraically independent over $F\av$.
\end{lemma}
\begin{proof}
By induction on $n$, it suffices to treat the case when $n = 1$.
Let $b \in K \setminus \acl(F\av)$.
Since $W \subseteq K$ is large, then there exists
$b_1, b_2, b_3, b_4 \in W$ such that $(b_{1} - b_{2})/(b_{3} - b_{4}) = b$ and $b_3 \not = b_4$.
Therefore, at least one of the $b_{i}$ is not in $\acl(F\av)$.
\end{proof}

\begin{thm}\label{dense}
There exists $\K \models T$ which is countable and of infinite rank over~$F$.
For any such $\K$, the set $\mathbb G$ is a dense subset of $\DerK$.
\end{thm}
Thus, in a precise topological sense,  $\Generic$ is a \textbf{generic} set
(notice that $\Generic$ is $\tauFO$-closed in~$\DerK$).
\begin{proof}
We have seen that each $I_{Z}$ is open.
It suffices to prove the following claim.
\begin{claim}
For every $Z \in \mathcal L$, $I_{Z}$ is dense.
\end{claim}
Let $Z \subseteq K^{n} \times K^{n}$.
Let $B_{\av, \bv}$ be a nonempty basic open set.
We have to verify that $I_{Z} \cap B_{\av, \bv}$ is nonempty.
Let $\eps \in B_{\av, \bv}$: that is, $\eps\in \DerK$ and $\eps \av = \bv$.
Let $\eps_{0}$ be the restriction of $\eps$ to $\acl (F\av)$.
Let $\cv \in \Pi_{n}(Z)$ be algebraically independent over~$F\av$.
We can extend $\eps_{0}$ arbitrarily to $\cv$; in particular, there exists
$\delta \in \DerK$ such that $\delta$ extends $\eps_{0}$ and $\delta \cv \in Z$.
Thus, $\delta \in I_{Z} \cap B_{\av, \bv}$.
\end{proof}

The following theorem gives a ``topological'' criterion for when a differential
system has a solution in models of $\Tdg$.
\begin{thm}
Let $\tuple{\K, \eps} \models \Td$.
Assume that $\K$ countable and of infinite rank over $F$.
Let $\av \in K^{\ell}$.
Let $\DerK(\av, \eps)$ be the set of  derivations $\delta$ on $K$ extending $\eta$ and
such that $\eps$ and $\delta$ coincide on $\Jet_{\epsilon}(\av)$ where it is the Jet related to the derivation $\epsilon.$

Let $Z \subseteq K^{n} \times K^{n}$ be $L$-definable with parameters~$\av$.
Let $\delta$ be a derivation in $K$ such that $\tuple{\K, \delta} \models \Tdg$.
\Tfae:
\begin{enumerate}[series=generic]
\item $I_{Z}$ is dense in $\DerK(\av, \eps)$;
\item $I_Z$ is nonempty;
\item $I_{Z} \cap \mathbb G$ is nonempty;
\item\label{en:generic-4} $\delta \in I_{Z}$.
\end{enumerate}
\end{thm}
\begin{proof}
First of all it is easy to see that  $(1) \implies (2),$ and $(4) \implies (3)$ are obvious. Moreover (2) is equivalent to (3), since $I_{Z}$ by Lemma \ref{open} is open and $\mathbb G$ by Theorem \ref{dense} is dense in $\DerK$.

We prove first the case when $\av$ is empty (that is, $Z$ is
$L$-definable without parameters and $\DerK(\av, \eps) = \DerK$).  

In this case, we can add another equivalent formulation to \ref{en:generic-4}:
\begin{enumerate}[resume=generic]
\item $\mathbb G \subseteq I_{Z}$.
\end{enumerate}

Since $\Tdg$ is complete, and ``$\delta \in I_{Z}$'' can be expressed as a
first-order sentence (without parameters), we have that 
 $(4) \implies  (5),$ the converse is trivial so (4), (5) are equivalent.
Therefore, (1) is equivalent to (2).

Let us consider now the case when $\av$ is non-empty.
Let $F' \coloneqq F[\Jet_{\epsilon}(\av)]$ and let $\eta'$ be the restriction of $\eps$ to $F'$.
We denote by $L' \coloneqq L(F)$, and $T' \coloneqq T \cup \Diag(F')$.

We can consider the theory ${T'}^{\delta}_{g}$ of generic derivations on $K$
extending $\eta'$: notice that $\tuple{\K,\delta} \models {T'}^{\delta}_{g}$.
We can apply the previous proof to $T'$, since $Z$ is now $L'$-definable in $K$
without parameters (notice that we need to modify the definition of $\mathbb G$,
since we are restricting the space of derivations to those extending $\eta'$:
however, we already proved the equivalence between (2) and (3)). We conclude in this way the proof.
\end{proof}

Barbina and Zambella \cite{BZ} deal with a similar situation: however, we cannot use
their results, since to apply them to our setting we would need that $\K$ is
countable and  saturated.
Maybe there could be a common refinement if one could weaken their assumption to
$\K$ resplendent (since every countable consistent theory has a countable
resplendent model: see \cite{Hodges}).

\section{Conjectures and open problems}

In the article \cite{FT:dexp}, we posed several conjectures. We now conclude the paper with a list of additional open problems, remarks, and some further ideas.

\subsection{Definable types}
Let $\tuple{\K, \deltabar}\models \Tdsg$.
Given a type $p \in S_{\Ldb}^{n}(\K)$, let $\av$ be a realization of $p$; we
define $\tilde p \in S_{L}^{\Gamma_n}(\K)$ as the $L$-type of $\Jet(\av)$ over~$\K$.
\begin{open problem}
Is it true that $p$ is definable iff $\tilde p$ is definable?
We conjecture that it is true when $\Tdsg = \Tdbg$ and $T$ has NIP.
\end{open problem}

\subsection{Zariski closure}
Given $X \subseteq K^{n}$, denote by $X^{Zar}$ be the Zariski closure of $X$.
\begin{questions}[See \cite{FLL}]
1) Let $\Pa{X_{i}: i \in I}$ be an $L$-definable family of subsets of $K^{n}$.
Is $\Pa{X_{i}^{Zar}: i \in I}$ also $L$-definable?

2) Assume that 1) holds for $\K$.
Let $\tuple{\K, \deltabar} \models \Tdsg$.
Let $\Pa{X_{i}: i \in I}$ be an $\Ldelta$-definable family of subsets of $K^{n}$.
Is $\Pa{X_{i}^{Zar}: i \in I}$ also $\Ldelta$-definable?
\end{questions}

We have affirmative answers to both questions, which will be presented in a paper in preparation.
\subsection{Kolchin polynomial}
Let $\tuple{\monster, \deltabar}$ be a monster model of~$\Tdbg$.
Let $\av \in \monster^{n}$, $B \subseteq \monster$ such that 
$\deltabar B \subseteq B$.
There exists a polynomial $\omega_{\av \mid B}(t)$ such
that, for $n$ large enough,
$\rk(\Jet_n(\av) \mid B) = \omega_{\av \mid B}(n)$, where $\Jet_n(\av) = \set{\mu \av: \mu \in \Gamma \mbox{ and } \abs \mu \leq n}$ (see
\cite{Kolchin}).
The degree of the polynomial is at most $k$; denote by $\mu(\av)$ the leading
monomial of $\omega_{\av \mid K}$ (including its coefficient).
Let  $X \subseteq K^{n}$ be $\Ldb$-definable with parameters~$\bv$: define
\[\begin{aligned}
\mu(X) &\coloneqq \sup \Pa{\mu(\av \mid {\Jet}(\bv)): \av \in X}\\
\omega_{X} &\coloneqq \sup \Pa{\omega_{\av \mid {\Jet(\bv)}}: \av \in X}.
\end{aligned}
\]
where the supremum for $\omega$ is taken inside $\R[t]$ w.r.t. the order $p>q$ if $\lim_{t \to + \infty} p(t) > q$ (\cf \cite{FLL}),
while the supremum of $\mu$ is taken \wrt the order
$r t^d \leq r' t^{d'}$ if $d < d'$ or $d = d'$ and $r \leq r'$
(\cf \cite{Fornasiero:Hilbert-length}).
Notice that, by \cref{sec:ind}, $\mu(X)$ and $\omega_{X}$ are well-defined (that is, they not
depend on the choice of the parameters~$\bv$).
\begin{conjecture}[See \cites{FLL, Riviere:09}]\label{conj:Kolchin}
The suprema in the definitions of $\omega_X$ and $\mu(X)$ are maxima.
Moreover, $\omega$ and $\mu$ are definable in families:
that is, for every $\Ldb$-definable  family $\Pa{X_{i}: i \in I}$ there exists
a partition of $I$ into finitely many definable set $I = I_{1} \sqcup \dotsb \sqcup I_{m}$
such that $\mu(X_{i})$ and $\omega_{X_{i}}$ are constant on each $I_{j}$.
\end{conjecture}

\begin{open problem}
What is the ``geometric'' meaning of $\mu(X)$?
Notice that, up to a multiplicative constant, the $k^{th}$ coefficent of $\omega_{X}$
is equal to $\deltadim(X)$.
\end{open problem}

If \cref{conj:Kolchin} is true, then the function $X \mapsto \mu(X)$ behaves
like a dimension on $\Ldelta$-definable sets (with the difference that the
values of $\mu$ are not natural numbers, but monomials: \cf \cite[\S2.5]{Dries:89}).
\begin{conjecture}
Assume that  $\monster$ is endowed with a topology $\tau$
satisfying some suitable conditions.
Let $\tau_{\deltabar}$ be the topology on $\monster$ induced by the embedding
$\monster \to \monster^{\Gamma}$, $x \mapsto  \Jet(x)$ (where $\monster^{\Gamma}$ is endowed with the
product topology induced by~$\tau$).
Denote by $\overline X^{\tau_{\deltabar}}$ 
the $\tau_{\deltabar}$-closure of~$X$.
Then, for every $X \subseteq \monster^{n}$ which is $\Ldb$-definable and nonempty,
$\overline X^{\tau_{\deltabar}}$ is also $\Ldb$-definable, and
$\mu(\overline X^{\tau_{\deltabar}}\setminus X) < \mu(X)$.
\end{conjecture}

\bigskip

\noindent \textbf{Acknowledgements} The authors thank Elliot Kaplan, Itay Kaplan, Noa Lavi,  Silvain Rideau-Kikuchi, and Marcus Tressl. We thank the anonymous referee for the careful reading of the manuscript
and for the constructive suggestions.

\bigskip

\noindent \textbf{Funding}  
Both  authors are members of the ``National Group for Algebraic and Geometric
Structures, and their Applications'' (GNSAGA-INDAM). 
The authors acknowledge financial support from INdAM (Istituto Nazionale di Alta Matematica), particularly through the Intensive Period in Model Theory, ``Model Theory and tame expansions of
topological Fields'', Naples, 19~May--18~July 2025.




\printbibliography


\end{document}